\documentclass[12pt]{article}
\usepackage{amsmath,amsfonts,amssymb}
\usepackage{graphicx,color}

\usepackage{theorem}

\theorembodyfont{\upshape}

\newcommand{\qed}{\nopagebreak\par\hspace*{\fill}$\square$}

\newcommand{\ZZ}{{\mathbb Z}}
\newcommand{\RR}{{\mathbb R}}

\newcommand{\CC}{{\mathbb C}}

\def\gg{{\mathfrak{g}}}
\def\FF{{\mathfrak{F}}}
\def\hh{{\mathfrak{h}}}

\numberwithin{equation}{section}

\title{Representations of Toroidal general linear Superalgebra}
\author{S.Eswara Rao\\
School of Mathematics\\
Tata Institute of Fundamental Research,\\
Mumbai, India.\\[2mm]
email: senapati@math.tifr.res.in}
\begin{document}

\maketitle

\begin{abstract}
We consider general linear superalgebra (type $A$) and tensor with 
Laurent polynomial ring in several variables. We then consider the 
universal central extension of this Lie superalgebra which we call 
toroidal superalgebra. We give a faithful representation of toroidal 
superalgebra using vertex operators and bosons.
\end{abstract}

\section*{Introduction}

The purpose of this paper is to construct a faithful representation of 
toroidal superalgebra. To define toroidal superalgebra, let $\gg$ be any 
finite dimensional Lie superalgebra with an invariant, supersymmetric and 
even bilinear form. Let $A = \CC [t_1^{\pm 1},\cdots, 
t_q^{\pm1}]$ be a Laurent polynomial ring in $q$  
commuting variables. Then $\gg \otimes A$ carries a natural Lie 
superalgebra structure. We now consider $\gg \otimes A \oplus 
\Omega_{A/d_A}$ where $\Omega_{A/d_A}$ is the module of Kahler 
differentials (see section 4 for details) which is central extension of 
$\gg \otimes A$ and is called toroidal superalgebras. They were first 
introduce in ${[EZ]}$ and ${[IK]}$. If we take $\gg$ to be a Lie algebra, 
we get the usual toroidal Lie algebra and they are extensively studied. 
See for example ${[E M] , [E M Y]}$, expositary article $[E]$ and the reference there in.\\
In this paper we specialize to the case $\gg = sl {(M|N)}$ (See 
section 2 for details). In this case the toroidal superalgebra is the 
universal central extention of $\gg \otimes A$ and we denote it by 
$\tau$. See ${[IK]}$ (We take $M > 1, N > 1$ and if $M = N$ we take $M > 
2).$\\

In this paper we give a representations for $\tau$ using vertex 
operators and bosons. For a suitable non-degenerate integral lattice 
$\overline{\Gamma}$ and a sublattice $Q$, we construct a Fock space  $V 
{[\overline{\Gamma}]}.$ For each $\alpha \in Q,$ define vertex operator 
$X {(\alpha ,z)}$ such that its Fourier components $X_n (\alpha)$ act 
on the Fock space. We also consider the Fock space $\FF$ given  by 
bosons (see ${[W1]})$. Then the product of vertex operators and bosons 
( they commute ) acting on the bigger Fock space $V {[\overline{\Gamma}] 
\otimes \FF}$ give a representation for the toroidal superalgebra 
$\tau$. It should be mentioned that, for the case $q =1,$ the 
representation is stated in ${[KW]}$ without any details whatsoever. 
Actually one needs a cocyle on the lattice $\overline{\Gamma}$ to get 
representation even in the case $q = 1$. We have worked out all the 
details in section 3.\

In section 1 we recall the standard construction of vertex operators
${[FK]}$ (homogeneous picture) and state the representation for the 
Lie algebra case. In section 2 we give a Chevalley type basis for
$ gl {(M|N)}$ with the help of a co-cycle $(2.4)$. In section $3$, we 
define bosons and the Fock space $\FF$ on which bosons act. We now consider the 
bigger Fock space $V {[\overline{\Gamma}]\otimes \FF}$ and make the 
affine superalgebra ( the toroidal superalgebra for $q = 1$ ) act on it 
(Proposition 3.3). In the last section we state and prove the main 
Theorem $4.6$ which is to give representation to the toroidal superalgebra $\tau$.

\section{Vertex operators and Fock space}
In this section we recall the basic construction $[FK]$ of vertex 
operators acting on Fock space. Unlike in  $[FK]$ our integral lattice 
need not be even. The lattice contains vectors of odd norm. See the 
books $[K 1]$ and $[Xu]$. All our vector space are over complex 
numbers $\CC$.

Let $d$ be a positive integer and let $\hh$ be a vector space with basis 
$e_1, \cdots, e_d.$ Let $( , )$ be a non-degenerate bilinear form on $\hh$ 
such that $ (e_i, e_j ) = \delta_{ij}.$ 
Let $\Gamma = \oplus^d_{i=1} \ZZ \ e_i$ be an integral lattice so that $\hh = 
\Gamma \displaystyle \bigotimes_\ZZ \CC$.\

Define a Heisenberg algebra.
$$
\widehat{\hh} = \displaystyle {\oplus_{k \in \ZZ}} \hh (k) \oplus \CC K
$$

Where each $\hh (k)$ is an isomorphic copy of $\hh$ and the isomorphism  
given by 
$\alpha 
\longmapsto \alpha (k).$\
The Lie algebra structure on $\widehat{\hh}$ is define by

$$
[ \alpha (k), \beta (m)] = k (\alpha, \beta) \delta_{k+m,0} 
K\leqno{{\bf(1.1)}}
$$
and $K$ is central

Define $\hh_\pm = \oplus_{k \gtrless 0} \hh(k)$\

The Fock space representation of $\hh_+ \oplus \hh_-$ is the symmetric 
algebra $S(\hh_-)$ of $\hh_-$ together with the action of $\hh_+ \oplus  
\hh_-$ on $S(\hh_-)$ defined by\\
$K$ acts as 1.\\
$\alpha (-m)$ acts as multiplication by $\alpha (-m)$ for $m >0$ \\
$\alpha(m)$ acts as the unique derivation on $S(\hh_-)$ for which $\beta 
(-n) \longmapsto \delta_{m,n} \ m (\alpha , \beta)$ for $m,n > 0$\\
Note that  $S(\hh_-)$ affords an irreducible representation of $\hh_+ \oplus 
\hh_-$.\\
Let $ \Gamma_{\overline{i}} = \{\alpha \in \Gamma | (\alpha, \alpha) \in 2 \ZZ + i\}$\\
Then clearly $\Gamma = \Gamma_{\overline{0}} \cup \Gamma_{\overline{1}}$\\
{\bf{(1.2)}} We now define a Co-cycle\\
$F : \Gamma \times \Gamma \longrightarrow \{\pm 1\}$ by\\
$F (0, \alpha) = F (\alpha ,0) =1$ \ and\\
$F (e_i, e_j) = \begin{cases} 1& if \ i\leq j\\
-1& if \ i> j \end{cases}$\\
We extend the map bimultiplicatively to $\Gamma$.\\
Note that $F (e_i,-e_j) = F (e_i,e_j)$
and $F (\alpha , \beta) = F (\alpha , \beta)^{-1}$ \\
We also note the following properties of $F$ which can be easily checked.\\
$F (\alpha , \beta) F (\alpha +\beta,\gamma) = F (\beta, \gamma) F 
(\alpha ,\beta + \gamma)$ for all $\alpha, \beta, \gamma \in \Gamma$\\
$F (\alpha, \beta) F (\beta, \alpha)^{-1} = (-1)^{(\alpha, \beta) + ij}.$
Where $\alpha \in \Gamma_{\overline{i}}$ and $\beta \in 
\Gamma_{\overline{j}}$.\\
To see this let $\alpha = \sum m_k e_k, \beta =\sum n_l e_l$\\
Then
\begin{eqnarray*}
\sum m_k^2 &= &(\alpha, \alpha) \equiv i (2)\\
\sum n_l^2& = &(\beta, \beta) \equiv j (2)
\end{eqnarray*}

and $(\alpha, \alpha) (\beta, \beta ) \equiv i j (2)$\\
Now consider\\
$F (\alpha, \beta) F (\beta, \alpha)^{-1} = F (\alpha, \beta) F (\beta, 
\alpha)$\\
$= \displaystyle{\prod_{k > l}} F (e_k, e_l)^{m_k n_l} 
\displaystyle{\prod_{l > k}} F (e_l, e_k)^{n_l m_k}$\\
$= (-1)^D$ where $D = \displaystyle{\sum_{k \not= l}} m_k n_l$\\
But $ \displaystyle{\sum_{k \not= l}} m_k n_l +  \displaystyle{\sum_{k} 
m_k n_k} =  \displaystyle{\sum_{k, l}} m_k n_l \equiv
\displaystyle{\sum_{k,l}} m_k^2 n_l^2 (2)$\\
So 
\begin{eqnarray*}
D &\equiv &-(\alpha , \beta) + (\alpha , \alpha) (\beta, \beta)\\
& \equiv& \big( (\alpha , \beta) + ij \big) (2)
\end{eqnarray*}
It will now follow
$F(\alpha , \beta) F (\beta, \alpha)^{-1} = (-1)^{(\alpha , \beta)+ ij}$\\
{\bf{(1.3)}} For each $\gamma$ in $\Gamma$, let $e^\gamma$ be a symbol and form the 
vector space $\CC[\Gamma]$ with basis  $\{e^\gamma,\gamma\in\Gamma\}$ over
 $\CC$. We define a twisted group algebra structure on $\CC[\Gamma]$ by

$$
e^\alpha \cdot e^\beta =F(\alpha, \beta) e^{\alpha + \beta}
$$
Consider the space $V[\Gamma] =\CC[\Gamma] \otimes S(\hh_{-})$\\
and define $\widehat{\hh}$ action on $V[\Gamma]$ by

\begin{eqnarray*}
\alpha(m) \cdot e^\gamma \otimes u & =& e^\gamma \otimes \alpha(m)u, 
\ m\neq 0\\
\alpha(0) \cdot e^\gamma \otimes u &=& (\alpha, \gamma) e^\gamma \otimes u
\end{eqnarray*}
It is a standard fact that $V[\Gamma]$ is $\widehat{\hh}$ -module.\\
$V[\Gamma]$ has a natural $\ZZ_2$-gradation.\\
A vector $e^\alpha\otimes u \in V[\Gamma]$ is called even if 
$(\alpha, \alpha) \equiv 0(2)$  and odd if $(\alpha, \alpha)\equiv 1(2)$.  Then 
let $V_{\bar{o}}[\Gamma]$ be the linear span of even elements and 
$V_{\bar{1}}[\Gamma]$ be the linear span of odd elements.  Then we have 
$V[\Gamma] =V_{\bar{o}}[\Gamma]\oplus V_{\bar{1}}[\Gamma]$.\\
An operator $T$ on $V[\Gamma]$ is called even operator if
$$
T(V_{\bar{i}}[\Gamma])\subseteq V_{\bar{i}}[\Gamma]
$$
and is called odd operator if
$$
T(V_{\bar{i}}[\Gamma]) \subseteq V_{\overline{i+1}} [\Gamma]
$$
Vertex operators\\
Let $z$ be complex valued variable and let $\alpha\in \Gamma$.\\
Define $T\pm(\alpha,z) =-\displaystyle\sum_{n\in\ZZ_\pm} \frac{1}{n}
\alpha(n)z^{-n}$\\
Define operator $z^{\alpha(0)}, \alpha \in \Gamma$, by $z^{\alpha(0)} \cdot
e^\gamma \otimes u = z^{(\alpha,\gamma)} e^\gamma \otimes u.$\\
Define Vertex Operator\\
$ Y(\alpha,z) =e^\alpha z^{\alpha(0)}\  \mbox{exp} T_{-}
(\alpha, z) \ \mbox{exp} \ T_+(\alpha, z)$\\
For $\alpha \in \Gamma_{\bar{0}}$  define
$$
X(\alpha, z) =z^{\frac{(\alpha,\alpha)}{2}} Y(\alpha, z) \ \mbox{and}
$$
write
$$
X(\alpha, z) =\displaystyle\sum_{n\in \ZZ} X_n(\alpha) z^{-n}
$$
Then $X_n(\alpha)$ is a even operator on $V[\Gamma]$.\\
For $\alpha \in \Gamma_{\bar{1}}$ define
$$
X(\alpha, z) =Y(\alpha, z)\ \mbox{and}
$$
write
$$
X(\alpha, z) = \displaystyle\sum_{n\in \ZZ} X_{n -\frac{1}{2}} (\alpha)z^{-n}
$$
So that $X_{n-\frac{1}{2}} (\alpha)$ is an odd operator.  For $\alpha \in 
\Gamma
$, define
$$
\alpha(z) = \displaystyle\sum_{n\in \ZZ} \alpha(n) z^{-n-1}
$$
which is an even operator.\\
{\bf{(1.4)}} We will now introduce delta function and recall some standard facts
from Section 2 of [FLM].\\
Define delta function
$$
\delta(z) =\sum_{n\in \ZZ} z^n
$$

Then  we have the following lemma holds. See [FLM] for proof and definition.
Suppose $X(z,w) =\displaystyle \sum_{m,n \in \ZZ} X_{m,n}z^m w^n,$
define $D_z X(z,w) \ = \sum m X_{m,n} z^{m-1} w^n$

\paragraph*{Lemma (1.5)}

\begin{eqnarray*}
1)  X(z,w)\delta(z/w) &= &X(w,w)\delta(z/w)\\
2)  X(z,w)\delta(z,w) &=& X(z,z)\delta(z/w)\\
3)  X(z,w)D_z(\delta(z/w))&=& X(w,w) D_z(\delta(z/w))-(D_z X)(w,w)
\delta(z/w)\\
4)  X(z,w) D_w(\delta(z/w)) &=& X(z,z)D_w(\delta(z/w))-(D_w X)(z,z) \delta(z/w)
\end{eqnarray*}
Let $\Delta =\{ \alpha_{ij} =e_i -e_j, i\neq j\} \subseteq \Gamma$\\
Then $\Delta$ is a finite root system of type $A$. The following is a well 
known result. See [FK] or [EM].

\paragraph*{Proposition (1.6)} The operator $\alpha_{ij} (m), X_n 
(\alpha_{ij}), \ K=Id$ for $i \neq j, m, n \in \ZZ$, defines a representation 
of an affine Lie-algebra of type A on $V[\Gamma]$.
\paragraph*{(1.7)} We will now introduce normal ordering and super commutator
on operators $Y(\alpha, z)$.\\
Suppose $\alpha \in \Gamma_{\bar{i}}$ and $\beta \in \Gamma_{\bar{j}}$
$$
: Y_n(\alpha) Y_m(\beta) : =\begin{cases} Y_n(\alpha) Y_m(\beta) & \mbox{if}
 \ n \leq m\\
(-1)^{ij} Y_m(\beta) Y_n(\alpha) &\mbox{if} \ n > m
\end{cases}
$$
$$
[Y_n(\alpha), Y_m(\beta)] =Y_n(\alpha) Y_m(\beta)-(-1)^{ij} Y_m(\beta) Y_n(\alpha)
$$

We now prove the following

\paragraph*{Lemma (1.8):} Let $\alpha \in \Gamma_{\bar{i}}, \beta\in
\Gamma_{\bar{j}}$\\
1) $[Y(\alpha,z), Y(\beta, w) =0 \ if \ (\alpha, \beta)\geq 0$\\
2) $[Y(\alpha, z), Y(\beta, w)]= F(\alpha, \beta) Y(\alpha + \beta, z) z^{-1}
 \delta (z/w) \ \mbox{if} \ \ (\alpha,\beta)=-1$\\
3) $[\alpha(z), Y(\beta, w)] =(\alpha, \beta) z^{-1} Y(\beta, z) \delta(z/w)$

\paragraph*{Proof} We have $z^{\alpha(0)} e^\beta =z^{(\alpha, \beta)} 
e^\beta z^{\alpha(0)}$\\
Note that $[T_{-} (\alpha, z), T_{+} (\beta, w)] =\log (1-w/z)^{(\alpha,\beta)}$
which will follow from (6.1.57)) of $[Xu]$.\\
Let $T(\alpha,\beta, z,w)= e^{\alpha+\beta} z^{\alpha(0)} w^{\beta(0)}\ 
exp \ T_{-}(\alpha, z) exp \ T_{-}(\beta, w)\ exp \ T_{+} (\alpha, z) \ exp \
T_{+} (\beta, w)$\\
Consider\\
$e^\alpha z^{\alpha(0)} \ exp \ T_{-} (\alpha, z) \ exp \ T_{+} (\alpha, z)$.	
$e^{\beta} z^{\beta(0)} \ exp \ T_{-} (\beta, w)\ exp \ T_{+} (\beta, w)$\\
$=F(\alpha, \beta) T(\alpha, \beta, z,w) z^{(\alpha,\beta)} (1-w/z)^{(\alpha , \beta)}$\\
By symmetry we have
$$
Y(\alpha, z) Y(\beta, w) -(-1)^{ij} Y(\beta, w) Y(\alpha, z)
$$
$$
=T(\alpha, \beta, z, w) (F(\alpha, \beta) z^{(\alpha, \beta)}
(1-w/z)^{(\alpha,\beta)}$$
$$
-(-1)^{ij} F(\beta, \alpha) w^{(\beta, \alpha)} (1-z/w)^{(\alpha,\beta)})$$
Consider the term in brackets which is equal to 
$$
F(\alpha, \beta) \left( z^{(\alpha,\beta)} (1-w/z)^{(\alpha,\beta)}-(-1)^{ij}
F(\alpha,\beta)^{-1} F(\beta, \alpha) w^{(\alpha,\beta)} (1-z/w)^{(\alpha ,
\beta)}\right)$$
From (1.2) we have 
$$
(-1)^{ij} F(\alpha, \beta)^{-1} F(\beta,\alpha) =(-1)^{(\alpha,\beta)}
$$ 
So the term in the brackets equal to 
$$
F(\alpha,\beta)\left((z-w)^{(\alpha,\beta)} -(-1)^{(\alpha,\beta)} (w-z)^
{(\alpha,\beta)}\right)
$$ 
which is equal to zero if $(\alpha,\beta)\geq 0$.  Thus we have proved (1) of 
the Lemma. 
Suppose $(\alpha,\beta)= -1.$ Now the term in the bracket equal to 

\begin{eqnarray*}
&F(\alpha,\beta) \left(z^{-1} (1-w/z)^{-1} +w^{-1} (1-z/w)^{-1}\right)\\
&= w^{-1}F(\alpha,\beta)\left(\frac{w}{z}(1-w/z)^{-1} + (1-z/w)^{-1}\right)\\
&= F(\alpha, \beta) w^{-1} \delta(z/w)
\end{eqnarray*}
Then we have 
$$[Y(\alpha,z), Y(\beta,w)]=F(\alpha,\beta)T(\alpha,\beta,z,w).w^{-1} \delta(z/w)
$$
Now by Lemma 1.5(2)
\begin{eqnarray*}
&=& F(\alpha,\beta)T(\alpha,\beta,z,z)z^{-1} \delta(z/w)\\
&=&F(\alpha,\beta)Y(\alpha + \beta, z)z^{-1} \delta(z/w)
\end{eqnarray*}
This proves the second part of the Lemma. For the third part of the lemma
we refer to (6.1.56) of $[Xu]$\\
We need the following \\
{\bf{Corollary (1.9)}}
$$
\begin{array}{llll}
(1)& \mbox{For} \ j\neq k, \ [X_{m-\frac{1}{2}} (e_i), X_n (e_j -e_k)]=\delta_{ik} F(e_i, e_j -e_k) X_{m+n-\frac{1}{2}} (e_j)\\
(2)& [X_{m-\frac{1}{2}} (e_i), X_{n+\frac{1}{2}} (-e_j)] = \delta_{ij} F(e_i, -e_j) \delta_{m+n,o} K.\\
(3)& \mbox{Suppose} \ \alpha,\beta  \in \Gamma, \ [\alpha(m), X_n (\beta)] =(\alpha, \beta) X_{m+n} (\beta)
\end{array}
$$
We now recall the following known facts. Recall $\alpha_{ij} =e_i -e_j$\\
${\bf{(1.10)}} \ (1) \ [X_m (\alpha_{ij}), X_n(-\alpha_{ij})] =$
$F(\alpha_{ij}, -\alpha_{ij})(\alpha_{ij} (m+n) +\delta_{m+n,0}\  m K)$\\[1mm]
2) For $ i \neq j$\\
$$
\sum_{k\in\ZZ} : X_{k+\frac{1}{2}} (e_i)X_{n-k-\frac{1}{2}} (-e_j):=F(e_i, -e_j) X_n(e_i-e_j)
$$

$$
\sum_{k\in\ZZ} : X_{n- k -\frac{1}{2}} (-e_j) X_{k+ \frac{1}{2}} (e_i) :=F(-e_j, e_i) X_n(e_i-e_j)
$$\\
%%%%%%%%%%BY ANU%%%%%%%%%%%%%%%%%%
3)$: X (e_i,z) X (-e_i,z) : = e_i(z)$\\[2mm]
In other words\\
$\displaystyle{\sum_{k \in \ZZ}} : X_{k +\frac{1}{2}} (e_i) X_{n -k 
-\frac{1}{2}} (-e_i) : = e_i(n)$\\
1.10(1) is very standard. For example see $[KF]$ or $[EM]$.\\
1.10(2) Follows from the definition of vertex operators.\\
Note that equations are equivalent as the operator at $LHS$ anti-commute 
and $F(e_i,-e_j) = -F(-e_j,e_i).$\\
1.10(3) Which is non-trivial can be found in (5.37) of $[FZ]$\\
{\bf{Proposition (1.11)}}\\
The operators $e_k(n), X_m(\alpha_{ij}),1 \leq i,j,k \leq d, \ n,m \in 
\ZZ, \ K=Id$ defines a representation of $\stackrel{\wedge}{gl_d}$ on 
the space $V[\Gamma]$. Where $\stackrel{\wedge}{gl_d}$ is the standard 
affinization of $gl_d$.\\
{\bf{Proof}} In view of Proposition 1.6, the only thing we need to check 
is\\
$$
[e_k(r), X_m(\alpha_{ij})]= (e_k,\alpha_{ij}) X_{m+r} (\alpha_{ij})
$$
This will follow from Corollary 1.9(3).
\section{Affine Superalgebras}
{\bf{(2.1)}} \ A lie superalgebra is a $\ZZ_2-$ graded vector space $\gg = 
\gg_{\overline{0}} \oplus \gg_{\overline{1}}$ equipped with $\CC$ 
bilinear map $[ , ] : \gg \times \gg \rightarrow \gg,$ called the Lie 
superbracket, satisfying the conditions.\\
$(1) \ [\gg_{\overline{i}} ,\gg_{\overline{j}}] \subseteq  \ \gg_{\overline{i+j}}$\\
$(2) \ [X,Y] = - (-1)^{ij} [Y,X]$\\
$(3) \ \big[[X,Y],Z\big] = \big[X, [Y,Z]\big] -(-1)^{ij} \big[Y, [X,Z]\big]$\\
For all homogeneous elements $X \in \gg_{\overline{i}}, Y \in \gg_{\overline{j}}$ and $Z \in \gg_{\overline{k}}$\\
The subspace $\gg_{\overline{o}}$ is called even and the subspace $\gg_{\overline{1}}$ is called odd.
It is easy to see that $\gg_{\overline{o}}$ is the usual Lie algebra and $\gg_{\overline{1}}$ is $\gg_{\overline{o}}$
module. The identity $(3)$ is called super Jacobi identity.\\
{\bf{(2.2)}} \ An important example of a Lie superalgebra is the space of all endomorphisms of a $\ZZ_2 -$ graded vector space
$V=V_{\overline{o}} \oplus V_{\overline{1}}$. We assume $\dim V_{\overline{o}} = M$ and  $\dim V_{\overline{1}} = N$.
The $\ZZ_2 -$ gradation on $V$ naturally induces a $\ZZ_2 -$ gradation on\\
$$
End (V) = (End(V))_{\overline{o}} \oplus (End(V))_{\overline{1}}  
$$
by letting
$$
End (V)_j = \{ f \in End (V) : f(V_k) \subseteq V_{k+j} \mbox{for \ all}\ k \in \ZZ_2 \}
$$
End $V$ becomes a Lie superalgebra with Lie super bracket\\
$$
[f,g]= f \circ g -(-1)^{ij} g \circ f
$$
for all $f \in End(V)_i, g \in End(V)_j$. We fix a basis of $V_{\overline{o}}$ say 
$e_1, \cdots e_M$ and a basis $e_{M+1}, \cdots e_{M+N}$ of $V_{\overline{1}}$. We denote the Lie superalgebra $End(V)$ by 
$gl (M|N)$ which is of type $A$.\\
An elements $f$ in $End(V)$ can be represented by a matrix form with respect to the above basis.\\
$$
f=
\begin{pmatrix}
A & B\\
C & D
\end{pmatrix}
$$
where $A$ is $M \times M$ matrix, $B$ is $N \times M$ matrix, $C$ is $M \times N$ matrix and $D$ is $N \times N$ matrix.\\
It is easy to see that\\
$$
\begin{pmatrix}
A & O\\
O & D
\end{pmatrix} \in End (V)_{\overline{o}} \ and
$$

$$
\begin{pmatrix}
O & B\\
C & O
\end{pmatrix} \in End (V)_{\overline{1}} \ 
$$
If $X \in End(V) = gl(M|N)$ is a homogeneous of deg $i$, then we denote $|X| =
 i$.\\
{\bf{(2.3)}} Define super trace of $f , str f= Trace \ A-Trace \ D.$\\
Let $< , >$ be a form on $End V$ by $< f, g> = str (f \circ g).$\\
The form $< , >$ satisfy by the following properties\\
$ (1) \ < X,Y > =(-1)^{|X| |Y|} < Y,X >$ (Super-symmetric)\\
$ (2) \ < [X,Y],Z > = < X, [Y,Z]>$ (invariant)\\
$ (3) \ < gl(M|N)_{\overline{o}} \ | \ \gg l (M|N)_{\overline{1}}> =0$ (even)\\
The root system of $gl(M|N)$ is \\
$$
\Delta = \{\alpha_{ij} = e_i-e_j, i \neq j, 1 \leq i, j \leq M+N \}
$$
{\bf{(2.4)}} We will now give a Chevalley type basis with the help of a Co-cycle for 
$gl(M|N).$\

Let $T$ be a vector space spanned by vector $T_{ij} , 1 \leq i, j \leq M+N$ We
call $T_{ij}$ even if $1 \leq i, j \leq M$ or $M+1 \leq i, j \leq M+N$. We call $T_{ij}$ \
odd if $1 \leq i \leq M$ and $j>M$ or $1 \leq j \leq M$ and $i>M$. Let $T_0$ be the linear span of even
elements and $T_1$ be the linear span of odd elements. Then we have $T= T_0 \oplus T_1$ is a 
$\ZZ_2-$ graded space. We would like to give Lie superalgebra structure on $T.$\\
Let $\Gamma = \oplus^M_{i=1} \ZZ e_i$ and let\\
$F : \Gamma \times \Gamma \rightarrow \{ \pm 1\}$ be the Co-cycle as defined in $(1.2)$.\\
Let $\alpha_{ij} = e_i - e_j.$ Note that $(\alpha_{ij}, \alpha_{kl}) \geq -2.$\\	
We will now define a super bracket on $T$.\\
$(T1) \ 1\leq i, j, k, l \leq M$\\
$i\neq j, k\neq l$\\
$[T_{ij}, T_{kl}]= \begin{cases} 0 \ \mbox{if} \ (\alpha_{ij}, \alpha_{kl}) \geq 0\\
F(\alpha_{ij}, \alpha_{kl}) T_{il} \ \mbox{if} \ j=k, l\neq i, (\alpha_{ij}, \alpha_{kl})=-1\\
F(\alpha_{ij}, \alpha_{kl}) T_{kj} \ \mbox{if} \ l=i, j\neq k, (\alpha_{ij}, \alpha_{kl})=-1\\
F(\alpha_{ij}, \alpha_{kl}) (T_{ii} -T_{jj})\ \mbox{if} \ l=i, j=k, (\alpha_{ij}, \alpha_{kl})=-2\\
\end{cases}$\\[3mm]
$i=j, k\neq l$\\
$[T_{ii}, T_{kl}] = (e_{i}, \alpha_{kl}) T_{kl}$\\
$i=j, k=l$\\
$[T_{ii}, T_{kk}] = 0$\\[3mm]
$(T2) \ 1\leq i, j, k, l \leq N$\\
$[T_{i+M,j+M} , T_{k+M,l+M}]= \delta_{jk} T_{i+M,l+M}-\delta_{il} \ T_{k+M,j+M}$\\[3mm]
$(T3) \ 1 \leq i \neq j, l \leq M, 1 \leq k \leq N$\\
$[T_{ij} , T_{k+M,l}] = \delta_{il} F(e_j,e_i) T_{k+M,j}$\\
$[T_{k+M,l} T_{ij}] = \delta_{il} F(e_i,e_j) T_{k+M,j}$\\
$[T_{ii} , T_{k+M,l}]= - \delta_{il} T_{k+M,l}$\\ 
$[T_{k+M,l} , T_{ii}] = \delta_{il} T_{k+M,l}$\\[3mm]
$(T4) \ 1\leq i, j, k \leq M, 1 \leq l \leq N$\\ 
$[T_{ij} , T_{k,l+M}] = \delta_{jk} F(e_i,e_j) T_{i,l+M}$\\
$[T_{k,l+M} , T_{ij}] = -\delta_{jk} F(e_i,e_j) T_{i,l+M}$\\[3mm]
$(T5) \ 1\leq i, j, k \leq N, 1 \leq l \leq M$\\ 
$[T_{i+M,j+M} , T_{k+M,l}] = \delta_{jk} T_{i+M,l}$\\
$[T_{k+M,l} , T_{i+M,j+M}] = -\delta_{jk} T_{i+M,l}$\\[3mm]
$(T6) \ 1\leq i, j, l \leq N, 1 \leq k \leq M$\\ 
$[T_{i+M,j+M} , T_{k,l+M}] = -\delta_{li} T_{k,j+M}$\\
$[T_{k,l+M} , T_{i+M,j+M}] = \delta_{li} T_{k,j+M}$\\[3mm]
$(T7) \ 1\leq i,l \leq N, 1 \leq j, k \leq M$\\ 
$[T_{i+M,j} , T_{k,l+M}] = \delta_{jk} T_{i+M,l+M} + F(e_k,e_j) \delta_{li} T_{kj}$\\
$[T_{k,l+M} , T_{i+M,j}] = \delta_{jk} T_{i+M,l+M} + F(e_k,e_j) \delta_{li} T_{kj}$\\[3mm]
$(T8) \ 1\leq i,j \leq M, 1 \leq k, l \leq N$\\ 
$[T_{ij} , T_{k+M,l+M}] = [T_{k+M,l+M} ,T_{ij}] = 0$\\[3mm]
$(T9) \ 1\leq i,k \leq N, 1 \leq j, l \leq M$\\ 
$[T_{i+M,j} , T_{k+M,l}] = [T_{k+M,l} ,T_{i+M,j}] = 0$\\[3mm]
$(T10) \ 1\leq i,k \leq M, 1 \leq j, l \leq N$\\ 
$[T_{i,j+M} , T_{k,l+M}] = [T_{k,l+M} ,T_{i,j+M}] = 0$

We claim that $T$ is a Lie superalgebra with the above super bracket. It is sufficient to check the super Jacobi identity.
This can be directly checked case by case but very tedious.\

Notice that the super bracket is antisymmetric for vectors of type even-even or even-odd. The super bracket is symmetric
for odd-odd vectors. This is consistent with Lie superalgebra definition.\

Let $\hh$ be the subspace spanned by $T_{ii} , 1 \leq i \leq M+N$. Clearly $\hh$ is abelian and $T$ decomposes with respect
to $\hh$.\\
So by comparing the root system of $T$ with that of $gl (M|N)$ we conclude that $T \cong gl (M|N)$. See $[K2]$ for details.\\
{\bf{(2.5)}} We now define bilinear form $(,)_T$ on $T$\\
$(1) \ 1\leq i,j,k,l \leq M$\\ 
$(T_{ij} , T_{kl})_T = F (\alpha_{ij},\alpha_{kl}) \delta_{jk} \delta_{li}$\\[3mm]
$(2) \ 1\leq j,k \leq M , 1\leq i,l \leq N$\\ 
$(T_{i+M,j} , T_{k,l+M})_T = -\delta_{jk} \delta_{li}$\\
$(T_{k,l+M} , T_{i+M,j})_T = \delta_{jk} \delta_{li}$\\[3mm]
$(3) \ 1\leq i,j,k,l \leq N$\\ 
$(T_{i+M,j+M} , T_{k+M,l+M})_T = -\delta_{jk} \delta_{li}$\\[3mm]
$(4)$ \ All other brackets are zero.\

Extend$(,)_T$ to whole of $T$ bilinearly. It is straight forward to check
that $(,)_T$ is a supersymmetric,invariant and even.\\
We will now define affine superalgebra
$$
\widehat{T} = T \otimes \CC [t, t^{-1}] \oplus \CC K
$$
Let $X(m) = X \otimes t^m , X \in T, m \in \ZZ.$\\
$[X(m), Y(n)] = [X,Y] (m+n) + (X,Y)_{T} \ \delta_{m+n,o} \ mK$\\
We will now write down the affine bracket explicitly on the genarators of $T$.\\
$(R1) \ 1\leq i,j,k,l \leq M$\\ 
$i\neq j, k\neq l$\\
$[T_{ij} (m), T_{kl} (n)]=\begin{cases} 0 \ \mbox{if} \ (\alpha_{ij}, \alpha_{kl}) \geq 0\\
F(\alpha_{ij}, \alpha_{kl}) T_{il} (m+n)\ \mbox{if} \ j=k, l\neq i, (\alpha_{ij}, \alpha_{kl})=-1\\
F(\alpha_{ij}, \alpha_{kl}) T_{kj} (m+n)\ \mbox{if} \ l=i, j\neq k, (\alpha_{ij}, \alpha_{kl})=-1\\
F(\alpha_{ij}, \alpha_{kl}) ((T_{ii} -T_{jj}) (m+n) + m \delta_{m+n,0} K) \mbox{if} \ l=i, j=k, (\alpha_{ij}, \alpha_{kl})=-2\\
\end{cases}$\\[3mm]
$i=j, k\neq l$\\
$[T_{ii}(m), T_{kl} (n)] = (e_{i}, \alpha_{kl}) T_{kl} (m+n)$\\
$i=j, k=l$\\
$[T_{ii} (m), T_{kk} (n)] = m \ \delta_{ik} \ \delta_{m+n,0} K$\\[3mm]
$(R2) \ 1\leq i,j,k,l \leq N$\\ 
$[T_{i+M,j+M} (m), T_{k+M,l+M} (n)] = $\\
$ \delta_{jk} T_{i+M,l+M} (m+n)\ -\delta_{li} T_{k+M,l+M} (m+n)\ - \delta_{jk} \delta_{il}\ m \ \delta_{m+n,o} K$\\[3mm]
$(R3) \ 1\leq i \neq j,l \leq M, 1 \leq k \leq N$\\ 
$[T_{ij} (m) , T_{k+M,l} (n)] = \delta_{il} F (e_j, e_i) T_{k+M,j} (m+n)$\\
$[T_{k+M,l}(n),T_{ij}(m)] = \delta_{il} F (e_i, e_j) T_{k+M,j} (m+n)$\\
$[T_{ii}(m),T_{k+M,l}(n)]= -\delta_{il} T_{k+M,l} (m+n)$\\
$[T_{k+M,l} (n),T_{ii}(m)]=\delta_{il} T_{k+M,l} (m+n)$\\[3mm]
$(R4) \ 1\leq i,j,k \leq M, 1 \leq l \leq N$\\ 
$[T_{ij} (m) , T_{k,l+M} (n)] = \delta_{jk} F (e_i, e_j) T_{i,l+M}(m+n)$\\
$[T_{k,l+M}(n),T_{ij}(m)] = -\delta_{jk} F (e_i, e_j) T_{i,l+M}(m+n)$\\[3mm]
$(R5) \ 1\leq i,j,k \leq N, 1 \leq l \leq M$\\ 
$[T_{i+M,j+M} (m), T_{k+M,l}(n) = \delta_{jk} T_{i+M,l} (m+n)$\\
$[T_{k+M,l}(n) , T_{i+M,j+M} (m)] = -\delta_{jk} T_{i+M,l} (m+n)$\\[3mm]
$(R6) \ 1\leq i,j,l \leq N, 1 \leq k \leq M$\\ 
$[T_{i+M,j+M} (m), T_{k,l+M}(n) = -\delta_{li} T_{k,j+M} (m+n)$\\
$[T_{k,l+M}(n) , T_{i+M,j+M} (m)] = \delta_{li} T_{k,j+M} (m+n)$\\[3mm]
$(R7) \ 1\leq i, l\leq N, 1\leq j, k\leq M$\\
$[T_{i+M,j} (m), T_{k,l+M} (n)] =$\\
$\delta_{jk} T_{i+M,l+M} (m+n) + F(e_k,e_j) \delta_{li} \ T_{kj} (m+n)$\\
$- \delta_{jk} \ \delta_{li} \ m \ \delta_{m+n,o} \ K\\$
$[T_{k,l+M} (n) , T_{i+M,j} (m)]=$\\
$\delta_{jk} \ T_{i+M, l+M} \ (m+n) + F(e_k,e_j) \ \delta_{li} \ T_{kj} \ (m+n) $\\
$+ \delta_{jk} \ \delta_{li} \ n \ \ \delta_{m+n,o} \ K.$\\[3mm]
$(R8) \ 1\leq i, j\leq M, 1\leq K, l\leq N$\\
$[T_{ij} (m), T_{k+M, l+M} (n)] = [T_{k+M, l+M} (n), T_{ij} (m)] = o$\\[3mm]
$(R9) \ 1\leq i, k\leq N, 1\leq j, l\leq M$\\
$[T_{i+M,j} (m), T_{k+M,l} (n)] = [T_{k+M,l} (n), T_{i+M,j} (m)] = o$\\[3mm]
$(R10) \ 1\leq i, k\leq M, 1\leq j, l\leq N$\\
$[T_{i,j+M} (m), T_{k,l+M} (n)] = [T_{k,l+M} (n) , T_{i,j+M} (m)] = o$
%%%%%%%%%%%%%%%%%%%%%%%%%%%%%%%%%%%%%%%%%%%%%%%ANU
\paragraph*{(2.6)}In this  subsection we will  state some identities on 
operators acting on a superalgebra.\\
A $\ZZ_2$-graded algebra is called super algebra. Let $A= A_0 \oplus A_1$
be $\ZZ_2$- graded algebra. The elements of $A_0$ are called even and the 
elements of $A_1$ are called odd. An operator on $A$ is called even operator
if it take even space to even space and  odd space to odd. An operator is called 
odd if it takes even space to odd and takes odd space to even.  An operator is
called homogeneous if it is even or odd. Suppose $X$ is  homogeneous 
operator then let
$$
\begin{array}{lll}
\mid X\mid =0 & if \ X \ is \ even\\
\mid X \mid = 1 & if \ \ X \ is \ odd
\end{array}
$$
Suppose $X$ and $Y$ are homogeneous operator on $A$.  \\
Define $[X,Y] =XY -(-1)^{\mid X\mid \ \mid Y\mid} YX$\\
Note that $[Y,X] =-(-1)^{\mid X \mid \ \mid Y \mid} [X,Y]$.

The following can easily be verified.  For example see [K2].\\
Suppose $X,Y,Z$ are homogeneous operators on $A$. Then \\[1mm]
\noindent
{\bf{(2.7)}}\ \ \ (1)$ [X,YZ] = [X,Y]Z + (-1)^{\mid X \mid \ \mid Y \mid} Y[X,Z]$\\
\hspace*{0.5in} (2)\ $[XY,Z] = X[Y,Z] +(-1)^{\mid Z \mid \ \mid Y \mid} [X,Z]Y$\\
We need the following Lemma.

\paragraph*{Lemma 2.8} Suppose $X_1$ and $X_2$ are odd operators and $Y_1$ and
$Y_2$ are even operators. Further suppose $[X_i, Y_j]=0$ for all $i$ and $j$.
Then
$$
[X_1 Y_1, X_2 Y_2] = [X_1, X_2] Y_1 Y_2- X_2 X_1 [Y_1, Y_2].
$$
\paragraph*{Proof} Director verification. Just expand the brackets.  Note that
$X_1 Y_1$ and $X_2 Y_2$ and odd operators.

\section{Bosons and Fockspace} In this section we recall the definitions of bosonic
fields acting on a Fock space $\FF$ from [W2].\\
For $1 \leq j \leq N$ let
$$
\begin{array}{lllll}
\varphi^j (z) &=&\sum_{r\in \ZZ} \varphi^j_{r-\frac{1}{2}} z^{-r} \ and \\
\varphi^{j^{*}} (z) &=&\sum_{r \in \ZZ} \varphi^{j^{*}}_{r-\frac{1}{2}} z^{-r}
\end{array}
$$
$2N$ bosonic fields and the coefficients  $\varphi^j_{r-\frac{1}{2}}, 
\varphi^{j^*}_{r-\frac{1}{2}}$ acting on a Fock space $\FF$ satisfying\\
{\bf{(3.1)}} \ \ (1) $\varphi^i_{r-\frac{1}{2}} \varphi^j_{s-\frac{1}{2}} -
 \varphi^j_{s-\frac{1}{2}} \varphi^i_{r-\frac{1}{2}}=0$\\
\hspace*{0.5in}(2)\ $\varphi^{i *}_{r-\frac{1}{2}}  \varphi^{j *}_{s-\frac{1}{2}} -
 \varphi^{j *}_{s-\frac{1}{2}}  \varphi^{i *}_{r-\frac{1}{2}} =0$\\
\hspace*{0.5in}(3)\ $\varphi^{i}_{r-\frac{1}{2}}  \varphi^{j *}_{s-\frac{1}{2}}-
  \varphi^{j*}_{s-\frac{1}{2}}  \varphi^{i}_{r-\frac{1}{2}} =
-\delta_{r+s-1,0} \delta_{ij}$\\
Then the fock space $\FF$ admits a vector $|0>$, called the vacuum
with the following properties.
$$
\varphi^i_{r-\frac{1}{2}} |0> =\varphi^{j^{*}}_{s-\frac{1}{2}} |0 > =0 \ \forall i, \
j, \forall \ r, s >0$$
Further $\varphi^i_{r-\frac{1}{2}}, \varphi^{j *}_{s-\frac{1}{2}}, \ r,s \leq 0$
act freely on $\FF$ and any element of the Fock space $\FF$ is linear
combinations of elements obtained by applying 
$$
\varphi^i_{r-\frac{1}{2}}, \varphi^{j *}_{s-\frac{1}{2}}, r,s \leq 0
$$
The operators $\varphi^i_{r-\frac{1}{2}}, \varphi^{j *}_{s-\frac{1}{2}}, r, s > 0$
(respectively $ r, s \leq 0$) are called annihilation (respectively creation) 
operators.

The purpose of this section is to give representation of $\stackrel{\wedge}{gl}(M/N)$
on the space $V[\Gamma] \otimes \FF$.  The  vertex operators $X_n(\alpha)
, \alpha(n)$ acting on the first component and the bosons $\varphi^i_{r-\frac
{1}{2}}$ and $\varphi^{j *}_{s-\frac{1}{2}}$ act on the second component. In 
particular the vertex operators and bosons commute.

Note that $V[\Gamma] \otimes \FF$ is a $\ZZ_2$ -graded space with $e^\alpha 
\otimes u \otimes v$ is even if $(\alpha,\alpha)$ is even and odd if
$(\alpha,\alpha)$ is odd.\\
{\bf{(3.2)}} We now give a representation of $\stackrel{\wedge}{T}$ on the space 
$V[\Gamma]\otimes \FF$ by the following rule
$$
\begin{array}{lllllll}
(1) & 1\leq i \leq M,\\ &T_{ii} (z) \mapsto e_i(z)
=: X(e_i,z) X(-e_i,z):.\\
(2) & i\neq j, \ 1 \leq i, j\leq M&\\
& T_{ij} (z) \mapsto X(e_i -e_j,z)\\
(3) & 1\leq i \leq M, \ 1\leq j\leq N &\\
&T_{i,j+M}(z) \mapsto X(e_i, z) \varphi^{j^*}(z)\\
(4) & 1\leq j\leq N, \ 1\leq i \leq M&\\
&T_{j+M,i}(z)\mapsto X(-e_i,z)\varphi^j(z)\\
(5) & 1\leq i, \ j\leq N &\\
&T_{i+M,j+M}(z) \mapsto :\varphi^{i} (z) \varphi^{j^*} (z) :\\
(6) &K \mapsto Id.\\ 
\end{array}
$$
Here normal ordering is the following 
$$
\begin{array}{lllll}
:\varphi^i_{r-\frac{1}{2}} \varphi^{j^*}_{s-\frac{1}{2}} : & = \varphi^i_{r-
\frac{1}{2}} \varphi^{j^*}_{s-\frac{1}{2}} \ if \ r\leq s\\
& = \varphi^{j^*}_{s-\frac{1}{2}} \varphi^i_{r-\frac{1}{2}}\ if \ r > s\\
\end{array}
$$
Note that normal ordering is not necessary in (3) and (4) as the operators
commute. Normal ordering for (1) is defined in the  earlier section.\\

We will now state the main result of this section.  It is stated in [KW] but
without any proof.  We could not find the details anywhere in the literature.
So we decided to offer the details of the proof. It should be mentioned that we need a Co-cycle
and this is not given in $[KW].$\\
For $ 1 \leq i \leq M, \ 1\leq j \leq N$\\
Write
$$
S_{i,j+M}(z) =X(e_i, z)\varphi^{j^*} (z)=\sum_{m\in\ZZ} S_{i,j+M} (m) z^{-m-1}
$$
So that
$$
\sum_{r\in \ZZ} X_{r-\frac{1}{2}} (e_i)\varphi^{j^*}_{m-r+\frac{1}{2}} =
S_{i,j+M} (m)
$$
For $1 \leq j \leq N, \ 1\leq i \leq M$\\
Write
$$
S_{j+M,i} (z) =X(-e_i, z) \varphi^j(z) =\sum_{m\in \ZZ} S_{j+M,i} (m) z^{-m-1}
$$
So that
$$\sum_{r\in\ZZ} X_{r-\frac{1}{2}} (-e_i) \varphi^j_{m-r+\frac{1}{2}} =S_{j+M,i}
(m)
$$
For $1 \leq i, \ j \leq N$\\
Write
$$
S_{i+M, j +M} (z) : = : \varphi^{i} (z) \varphi^{j^*} (z) : =\sum_{m\in\ZZ} S_{i+M, j +M} (m) z^{-m-1}
$$
So that
$$
\sum_{r\in\ZZ}: \varphi^i_{r-\frac{1}{2}} \varphi^{j^*}_{m-r+\frac{1}{2}}:=
S_{i+M,j+M}(m)
$$
\paragraph*{Proposition 3.3} Notation as above.  The above map defines a 
representation of $\stackrel{\wedge}{T} \cong \stackrel{\wedge}{gl}(M/N)$ on the space $V[\Gamma]
\otimes \FF$.

\paragraph*{Proof} The only thing that we need to check is relations $(R1)$ to
$(R10)$ in {\bf{(2.5)}}.\\
Let $\RR$ be the field of real numbers.\\
Let 
$$
\begin{array}{llll}
\psi : \RR &\rightarrow & \{ 0,1\} \ such \ that\\
\psi(x) &=& 1 \ if \ \mid x \mid \leq 1\\
&=& 0 \ if \ \mid x \mid  >1
\end{array}
$$
We use the principle of truncation to deal with infinite sums. 
Consider $S_{ij}(m)$ which is infinite sum but $S_{ij} (m)\psi (\epsilon m)$ is finite
and tends to $S_{ij}(m)$ as $\epsilon \rightarrow 0$.\\
The relation is (R1) is standard which is nothing but standard vertex operator
construction (homogeneous picture). For example see [FK] or [EM].  We will first
verify (R7) which is most important relation.  For that it is sufficient to 
check the following:
$$
\begin{array}{lllll}
3.3(1) & [S_{i+M,j} (m), S_{k,l+M} (n)]\\
&= \delta_{jk} S_{i+M,l+M} (m+n) +F(e_k, e_j) \delta_{li} X_{m+n} (e_k -e_j)\\
&-\delta_{jk} \delta_{li} \cdot m \delta_{m+n,0},\mbox{if} \ j\neq k\\
&= \delta_{jk} S_{i+M,l+M} (m+n) + \delta_{li}\ e_k(m+n)\\
&-\delta_{jk}\  \delta_{li}\  m \ \delta_{m+n,0} \ if \ j=k\\[2mm]
3.3(2) & [S_{k,l+M} \ (n),\  S_{i+M,j} (m)]\\
&=\delta_{jk} S_{i+M,l+M}(m+n) +F(e_k, e_j) \delta_{li} X_{m+n} (e_k-e_j)\\
&+ \delta_{jk} \delta_{li} \ n  \ \delta_{m+n,0} \ if \ j\neq k\\
&= \delta_{jk} S_{i+M,l+M} (m+n) + \delta_{li} \ e_k(m+n)\\
&+ \delta_{jk} \delta_{li} \ n \ \delta_{m+n,0} \ if \ j=k
\end{array}
$$
We will first verify (R7) which is most important.\\
First consider
$$A= [S_{i+M,j}(m), S_{k,l+M}(n)]\ \psi(\epsilon n) =
\sum_r \sum_s\left[X_{r-\frac{1}{2}} (-e_j)\varphi^i_{m-r+\frac{1}{2}}, 
X_{s-\frac{1}{2}} (e_k)\varphi^{l^*}_{n-s +\frac{1}{2}}\right] \psi(\epsilon n)
$$
Recall that $X$ operator are odd and $\varphi$ operators are even and $X$ 
operators commute with $\varphi$ operators. So we can use Lemma (2.8). So
$$
A=\sum_{r,s\in\ZZ} \left[X_{r-\frac{1}{2}}(-e_j), X_{s-\frac{1}{2}}(e_k)\right]
\varphi^i_{m-r+\frac{1}{2}}\varphi^{l^*}_{n-s+\frac{1}{2}} \psi(\epsilon n)\\
$$

$$
-\sum_{r,s\in\ZZ} X_{s-\frac{1}{2}}(e_k)X_{r-\frac{1}{2}} (-e_j) \left[\varphi^i_{m-r+
\frac{1}{2}}, \varphi^{l^*}_{n-s+\frac{1}{2}} \right]\psi (\epsilon n)
$$
We will now use 3.1(3) and Corollary 1.9(2) to conclude
$$
\begin{array}{llll}
A&=\delta_{jk} \sum_{r,s\in\ZZ} F(e_k, -e_j) \delta_{r+s-1,0} \ \varphi^i_{m-r+\frac{1}{2}} \varphi^{l^*}_{n-s+\frac{1}{2}} \psi (\epsilon n)\\
&+ \delta_{il} \sum_{r,s\in \ZZ} X_{s-\frac{1}{2}} (e_k)X_{r-\frac{1}{2}}
(-e_j) \delta_{m+n-(r+s)+1,0} \psi(\epsilon n)\\
&=\delta_{jk} F(e_k,-e_j) \sum_{s\in\ZZ} \ \varphi^{i}_{m+s-\frac{1}{2}} \ \varphi^{l^*}_{n-s+\frac{1}{2}} 
\psi(\epsilon n)\\
&+ \delta_{il} \sum_{s\in\ZZ} X_{s-\frac{1}{2}} (e_k) X_{m+n-s+\frac{1}{2}}
(-e_j) \psi(\epsilon n)
\end{array}
$$
We will now normal  order the sums. For this we split the first term according
to $m+s -\frac{1}{2} \leq n-s+\frac{1}{2}$ (which is in the normal ordering)
or $m+s-\frac{1}{2} > n-s+\frac{1}{2}$. Similarly we split the second sum 
according to $s-\frac{1}{2}\leq m+n-s+\frac{1}{2}$ or $s-\frac{1}{2} > m+n-
 s+\frac{1}{2}$\\
Thus
$$
\begin{array}{llll}
A&=& \delta_{jk} \displaystyle{\sum_{s \epsilon \ZZ}}F(e_k,-e_j) :\varphi_{m+s-\frac{1}{2}} \ \varphi^{l^*}_{n-s+\frac{1}{2}}:\psi(\epsilon n)\\
&-& \delta_{jk} \delta_{il} \sum_{2s> n-m+1} F(e_k, -e_j) \psi(\epsilon n) \delta_{m+n,0}\\
&+& \delta_{il} \displaystyle{\sum_{s \epsilon \ZZ}} : X_{s-\frac{1}{2}} (e_k) X_{m+n-s+\frac{1}{2}} (-e_j):\psi(\epsilon n)\\
&+& \delta_{il} \ \delta_{jk} \ \sum_{2s > m+m+1} \psi(\epsilon n) \delta_{m+n,0}\\[3mm]
\end{array}
$$

$$
\begin{array}{llll}
A&=& \delta_{jk} S_{i+M,l+M} \psi(\epsilon n)\\
&-& \delta_{jk} \delta_{il} \sum_{s> n+\frac{1}{2}} \psi(\epsilon n) \delta_{m+n,0}\\
&+& \delta_{il} \displaystyle{\sum_{s \in \ZZ}}:X_{s-\frac{1}{2}} (e_k) \ X_{m+n-s+\frac{1}{2}} (-e_j):  \psi(\epsilon n)\\
&+& \delta_{il} \delta_{jk} \sum_{s>\frac{1}{2}} \psi(\epsilon n) \delta_{m+n,0}
\end{array}
$$
Now it is easy to see that
$$
\begin{array}{llll}
&\sum_{s>\frac{1}{2}} \psi(\epsilon n) \delta_{m+n,0} -\sum_{s>n+\frac{1}{2}} \psi(\epsilon n)
\delta_{m+n,0}\\
&=n \delta_{m+n,0} \psi(\epsilon n)
\end{array}
$$
Now taking $\epsilon  \rightarrow 0$ we see that for $j\neq k$  \ \ \ \mbox{(Using \ Corollary \ 1.9(1))}
$$
A\rightarrow  \delta_{jk} S_{i+M, l+M} (m+n) +\delta_{il} F(e_k, -e_j) X_{m+n} (e_k-e_j) + \delta_{jk} \ \delta_{li}\ n \delta_{m+n,0}
$$
Suppose $j=k$  \ \ \ \mbox{(Using \ 1.10(3))}
$$
\begin{array}{llll}
 A&=& \delta_{jk} S_{i+M, l+M} (m+n) \psi(\epsilon n)\\
&+& \delta_{il} \ e_j (m+n) \psi (\epsilon n)\\
&+& \delta_{jk} \ \delta_{li} \ n \delta_{m+n,0} \psi(\epsilon n)\\
\end{array}
$$
Now taking $ \epsilon \rightarrow 0$ we see that 
$$
\begin{array}{llll}
A & \rightarrow \delta_{jk} \ S_{i+M, l+M} (m+n) +\delta_{il} \ e_j(m+n) \\
& + \delta_{jk} \ \delta_{li} \ n \delta_{m+n,0}
\end{array}
$$
Thus we have verified 3.3(1).\\
We will now verify 3.3(2).\\
Consider
$$
\begin{array}{llll}
B &=& [S_{k,l+M} (n), S_{i+M,j} (m)] \psi (\epsilon n)\\
&=&\sum_{r\in\ZZ} \sum_{s\in\ZZ} \left[X_{s-\frac{1}{2}}(e_k)\varphi^{l^*}_{n-s+
\frac{1}{2}}, X_{r-\frac{1}{2}} (-e_j) \varphi^i_{m-r+\frac{1}{2}} \right]\psi(\epsilon n)\\
\end{array}
$$
We will now use Lemma (2.8).
$$
\begin{array}{llll}
B&=& \sum_{r,s\in\ZZ} \left[X_{s-\frac{1}{2}} (e_k), X_{r-\frac{1}{2}}(-e_j)\right] 
\varphi^{l^*}_{n-s+\frac{1}{2}}\varphi^i_{m-r+\frac{1}{2}} \psi(\epsilon n)\\
&-& \sum_{r,s\in\ZZ} X_{r-\frac{1}{2}}(-e_j) X_{s-\frac{1}{2}} (e_k) 
\left[\varphi^{l^*}_{n-s+\frac{1}{2}},\varphi^i_{m-r+\frac{1}{2}} \right]\psi(\epsilon n)\\
\end{array}
$$
We will now use 3.1(3) and Corollary 1.9(2)
$$
\begin{array}{llll}
B&=& \delta_{kj} \sum_{r,s\in\ZZ} \delta_{r+s-1,0} \varphi^{l^*}_{n-s+
\frac{1}{2}}\varphi^i_{m-r+\frac{1}{2}} \psi(\epsilon n)\\
&-& \delta_{li}\sum_{r,s \in\ \ZZ} X_{r-\frac{1}{2}} (-e_j) X_{s-\frac{1}{2}}
(e_k) \delta_{m+n-(r+s)+1,0} \psi(\epsilon n)\\
&=& \delta_{kj} \sum_{r\in\ZZ} \varphi^{l^*}_{n+r-\frac{1}{2}} \varphi^i_{m-r+\frac{1}{2}}\ \psi(\epsilon n)\\
&-&\delta_{li} \sum_{r\in\ZZ} X_{r-\frac{1}{2}} (-e_j) X_{m+n-r+\frac{1}{2}}(e_k)\psi(\epsilon n)
\end{array}
$$
We will now normal order the sums.  For that we split first sum according
$n+r-\frac{1}{2}\leq m-r+\frac{1}{2}$ or $n+r-\frac{1}{2} > m-r+\frac{1}{2}$
and split the second sum according $r-\frac{1}{2} \leq m+n -r+\frac{1}{2} \
or \ r-\frac{1}{2} > m+n-r+\frac{1}{2}$
$$
\begin{array}{llll}
B&=& \delta_{kj} \sum_{r\in\ZZ} :\varphi^{l^*}_{n+r-\frac{1}{2}} \ \varphi^i_{m-r+\frac{1}{2}}:
\psi(\epsilon n)\\
&+& \delta_{kj} \delta_{il} \sum_{2r>m-n+1} \psi(\epsilon n) \delta_{m+n,0}\\
&-& \delta_{il}\sum_{r\in\ZZ} : X_{r-\frac{1}{2}} (-e_j) X_{m+n-r+\frac{1}{2}}
(e_k) :\psi(\epsilon n)\\
&-& \delta_{il} \delta_{kj} \sum_{2r > m+n+1} \psi(\epsilon n) \delta_{m+n,0}
\end{array}
$$
The second + the last term equals
$$
\begin{array}{llll}
&\delta_{kj} \delta_{il} \left(\sum_{r>m+\frac{1}{2}} \psi(\epsilon n) \delta_{m+n,0}
-\sum_{r>\frac{1}{2}} \psi(\epsilon n) \delta_{m+n,0} \right)\\
& =-\delta_{kj}\delta_{il} m \delta_{m+n,0} \ \psi(\epsilon n)
\end{array}
$$
Note that
$$
\begin{array}{llll}
&-\sum_{r\in\ZZ} :X_{r-\frac{1}{2}} (-e_j) X_{m+n-r+\frac{1}{2}} (e_k):
&=\sum:X_{m+n-r+\frac{1}{2}}(e_k) X_{r-\frac{1}{2}}(-e_j):\\
&&=e_k(m+n) \ \mbox{if} \ j=k\\
&& =F(e_k,-e_j) X_{m+n} (e_k-e_j) \ if \ j\neq k
\end{array}
$$
Thus as $\epsilon \rightarrow 0$
$$
\begin{array}{llll}
B \rightarrow & \delta_{kj} S_{i+M, l+M}\  (m+n)+ \delta_{li} F(e_{k}, -e_j) X_{m+n}
(e_k-e_j) \ if \ k\neq j\\
B \rightarrow & \delta_{kj} S_{i+M, l+M} (m+n) +\delta_{li} e_k(m+n)\\
& -\delta_{kj} \delta_{li} m \delta_{m+n,0} \ \ if \ k=j
\end{array}
$$
This proves 3.3(2).\\
We will now verify the first part of (R3).  Consider for $ i\neq j$
$$
\begin{array}{llll}
[X_m(e_i-e_j), S_{k+M,l}(n)] =\sum_{r\in\ZZ} [X_m(e_i-e_j),X_{r-\frac{1}{2}}(-e_l) \varphi^k_{n-r+\frac{1}{2}}]\\
=\sum_{r\in\ZZ} [X_m(e_i-e_j), \ X_{r-\frac{1}{2}}(-e_l)] 
\varphi^k_{m-r+\frac{1}{2}} \ \ \ \mbox{(Using \ 2.1(3))}\\
=- \delta_{il} \sum_{r\in\ZZ} F(-e_l, e_i-e_j) X_{m+r-\frac{1}{2}} (-e_j) 
\varphi^{k}_{m-r+\frac{1}{2}} \ \ \ \mbox{(Using \ Corollary \ 1.9(1))}\\
=\delta_{il} F(e_j,e_i) S_{k+M,j} (m+n) $\qed$
\end{array}
$$
We will now verify the second part of (R3). For $ i\neq j$.
$$
\begin{array}{llll}
[S_{k+M,l} (n), X_{m}(e_i-e_j) =\sum_{r\in\ZZ}[X_{r-\frac{1}{2}}(-e_l)
\varphi^k_{n-r+\frac{1}{2}},\  X_m(e_i-e_j)]\\
=\sum_{r\in\ZZ}[X_{r-\frac{1}{2}} (-e_l), X_m (e_i-e_j)] \varphi^k_{n-r+\frac{1}{2}}\ \ \ \mbox{(Using \ 2.1(3))}\\
=\delta_{il} \sum_{r\in \ZZ} F (e_i,e_j)X_{m+r-\frac{1}{2}}(-e_j) \varphi^k_{n-r+
\frac{1}{2}} \ \ \ \mbox{(Using \ Corollary \ 1.9(1))}\\
=\delta_{il} F(e_i, e_j) S_{k+M,j} (m+n) \ (Here  \ and \ above \ we \ are\ 
using\\
the \ property \ of \ the\  co-cycle\ (1.2).)
\end{array}
$$
%&&&&&&&&&&&&&&&&&&&&anu&&&&&&&&&&&&&&&&&&&&&&&&&&&&&&&&&&&&&&&&&&&&&&&&%
We will now verify the third part of (R3).\\
$[e_i(m),S_{k+M,l}(n)]=\displaystyle {\sum_{\gamma \in \ZZ}} \left[e_i(m) ,X_{r -\frac{1}{2}} (-e_l) \varphi^k_{n-r+ \frac{1}{2}} \right]$\\
$=\displaystyle{\sum_{r \in \ZZ}} [e_i(m) ,X_{r -\frac{1}{2}} (-e_l)] \varphi^k_{n-r + \frac{1}{2}} \ \ \ \mbox{(Using \ 2.7)}$\\
$-\delta_{il} \displaystyle {\sum_{r \in \ZZ}} X_{r+m-\frac{1}{2}} (-e_l) \varphi^k_{n-r + \frac{1}{2}}\ \ \ \mbox{(Using \ Corollary \ 1.9(3))}$\\ 
$=-\delta_{il} \ S_{k+M,l}(m+n)$\\[2mm]
We will now prove the last part of (R3)\\
$[S_{k+M,l}(n), e_i(m)]=\displaystyle {\sum_{r \in \ZZ}}[X_{r- \frac{1}{2}} (-e_l) \varphi^k_{n-r + \frac{1}{2}},e_i(m)]$\\
$=\displaystyle {\sum_{r\in \ZZ}} [X_{r -\frac{1}{2}} (-e_l),e_i(m)]\varphi^k_{n-r + \frac{1}{2}} \ \ \ \mbox{(Using \ 2.7)}$\\
$=\delta_{il} \displaystyle {\sum_{r \in \ZZ}} X_{m+r -\frac{1}{2}} (-e_l) \varphi^k_{n-r + \frac{1}{2}} \ \ \ \mbox{(Using \ Corollary \ 1.9(3))}$\\
$=\delta_{il} \ S_{k+M,l}(m+n)$\\
The relations $R4, R5$ and $R6$ are similar and much easier as the central term does not appear.
We omit the details. The relations $R8, R9$ and $R10$ are trivial to verify.\\
We will now verify $R2$ which is actually very standard. We could not find any  decent proof which is suitable for our situation.
Any How we thought for the sake of completeness, we offer a short proof.\\
Consider
\begin{eqnarray*}
D&=& [S_{i+M,j+M} (m), S_{k+M, l+M} (n)] \psi (\epsilon n)\\
&=& \sum_{\gamma,s \in \ZZ} \left[\varphi^{i}_{r-\frac{1}{2}} \varphi^{j^*}_{m-r+ \frac{1}{2}}, 
\varphi^{k}_{s -\frac{1}{2}} \varphi^{l^*}_{n-s +\frac{1}{2}} \right] \psi (\epsilon n)
\end{eqnarray*}
We will be using $(2.7)$ and $(3.1)$. Recall that all $\varphi$ operators are even. Normal ordering is not necessary as they 
differ by a central operator.\\
\begin{eqnarray*}
D&=& \sum_{r,s \in \ZZ} \varphi^{i}_{r -\frac{1}{2}} \left[ \varphi^{j^*}_{m-r +\frac{1}{2}}, \varphi^{k}_{s -\frac{1}{2}} 
\varphi^{l^*}_{n-s +\frac{1}{2}} \right] \psi (\epsilon n)\\
&+& \sum_{r,s \in \ZZ} \left[\varphi^{i}_{r -\frac{1}{2}},  \varphi^{k}_{s -\frac{1}{2}} \varphi^{l^*}_{n-s +\frac{1}{2}} \right] 
\varphi^{j^*}_{m-r +\frac{1}{2}} \psi (\epsilon n)\\
&=&\sum_{r,s \in \ZZ} \varphi^{i}_{r -\frac{1}{2}} \left[ \varphi^{j^*}_{m-r +\frac{1}{2}}, \varphi^{k}_{s -\frac{1}{2}} \right]
\varphi^{l^*}_{n-s +\frac{1}{2}} \psi (\epsilon n)\\
&+& \sum_{r,s \in \ZZ} \varphi^{k}_{s -\frac{1}{2}} \left[\varphi^{i}_{r -\frac{1}{2}},\varphi^{l^*}_{n-s +\frac{1}{2}} \right]  
\varphi^{j^*}_{m-r +\frac{1}{2}} \psi (\epsilon n)
\end{eqnarray*}
\begin{eqnarray*}
D&=& \delta_{jk} \sum_{r,s \in \ZZ} \varphi^{i}_{r -\frac{1}{2}} \varphi^{l^*}_{n-s +\frac{1}{2}}
\delta_{m-r+s,0} \ \psi (\epsilon n)\\
&-& \delta_{il} \sum_{r,s \in \ZZ} \varphi^{k}_{s -\frac{1}{2}} \varphi^{j^*}_{m-r +\frac{1}{2}}
\delta_{n-s+r,0} \ \psi (\epsilon n) \ \mbox{(Using \ (3.1))}\\
&=& \delta_{jk} \sum_{r \in \ZZ} \varphi^{i}_{r -\frac{1}{2}} \varphi^{l^*}_{m+n-r +\frac{1}{2}}
\ \psi (\epsilon n)\\
&-&\delta_{il} \sum_{r \in \ZZ} \varphi^{k}_{n+r -\frac{1}{2}} \varphi^{j^*}_{m-r +\frac{1}{2}}
 \ \psi (\epsilon n)
\end{eqnarray*}
We need to normal order the sums. For this we need split the first sum accordingly to 
$r -\frac{1}{2} \leq n+m-r +\frac{1}{2}$ or $r -\frac{1}{2} > n+m-r +\frac{1}{2}$
We split the second sum according to $ n + r -\frac{1}{2} \leq m-r +\frac{1}{2}$ or $ n + r -\frac{1}{2} > m-r +\frac{1}{2}$\\
\begin{eqnarray*}
D&=& \delta_{jk} \sum_{r \in \ZZ} : \ \varphi^{i}_{r -\frac{1}{2}} \ \varphi^{l^*}_{n+m-r +\frac{1}{2}} :
\ \psi (\epsilon n)\\
&-&\delta_{jk} \delta_{ii} \sum_{2r> \ n+m+1} \psi(\epsilon n) \delta_{m+n,0}\\
&-& \delta_{il} \sum_{r \in \ZZ} : \ \varphi^{k}_{n+r -\frac{1}{2}} \ \varphi^{j^*}_{m-r +\frac{1}{2}} :
\ \psi (\epsilon n)\\
&+&\delta_{jk} \ \delta_{li} \sum_{2r> \ m-n +1} \psi(\epsilon n) \ \delta_{m+n,0}\\
\end{eqnarray*}
The second term and fourth term equal to $-\delta_{jk} \ \delta_{il} \ m\psi(\epsilon n)$\\
So as $\epsilon \rightarrow 0$
$$
D \rightarrow \delta_{jk} \ S_{i+M,l+M} (m+n) - \delta_{il} \ S_{k+M,j+M} (m+n) -\delta_{jk} \ \delta_{il} \ m \ \delta_{m+n,0.}
$$
This completes the Proof of the Proportion $3.3$
\section{Toroidal Superalgebras}
\quad In this section we define toroidal superalgebras corresponding to $gl(M|N)$ and give representation using a 
bigger Fock space.

Let $\gg$ be a Lie superalgebras and let $(\ ,\ )$ be a bilinear form on $\gg$ which supersymmetric, invariant and even.
Fix a positive integer $q$ and consider $A= \CC [t^{\pm 1}_{1} , \cdots t^{\pm 1}_{q}]$ a Laurent polynomial ring in $q$
commuting variables. Let $\overline{m}= (m_{1}, \cdots m_{q}) \in \ZZ^{q}$ and let 
$t^{\overline{m}}=t_{1}^{m_1} \cdots t_{q}^{m_q}$. Then $\gg \otimes A$ has a natural structure of Lie superalgebra.
Let $X(\overline{m}) = X \otimes t^{\overline{m}} \in \gg \otimes A$.

Let $ Z = \Omega_{A}/ d_{A}$ be the space of differentials so that $\Omega_{A}$ is spanned by symbols $t^{\overline{m}}
K_i, {\overline{m}} \in \ZZ^{q}, 1\leq i\leq q $ and $d_{A}$ is spanned by $\sum m_{i} t^{\overline{m}} K_{i}.$ We now
define Lie superalgebra structure on\	
$$
\gg \otimes A \oplus \Omega_{A}/ d_{A}
$$
$$
[X(\overline{m}), Y(\overline{n}) = [X,Y] (\overline{m} + \overline{n}) + (X,Y) d (t^{\overline{m}})t^{\overline{n}}
\leqno{\bf{(4.1)}} 
$$
$$
d(t ^{\overline{m}}) t^{\overline{n}} = \sum m_{i} \ t^{\overline{m}+ \overline{n}} K_i \leqno{\mbox{Where}} 
$$
and $X,Y \in \gg, \overline{m},\overline{n} \in \ZZ^{q}$.\\
$\Omega_{A}/ d_{A}$ is central.\\
We call this algebra a toroidal superalgebra. See $[E Z]$ for more details.\\
Define $sl (M|N) = \{ X \in g l (M|N) | str X = 0 \}$ which is a Lie supersubalgebra of $g l (M|N)$\\

\noindent
{\bf (4.2) Theorem ( Theorem (4.7) of [IK]).}\\
$sl (M|N) \otimes A \oplus \Omega_{A}/ d_{A}$ is the universal central extension of $sl (M|N) \otimes A$
(If $M=N,$ then we take $M>2$)

The purpose of this section is to give representation of $gl(M|N) \otimes A \oplus \Omega_{A}/ d_{A}$.
In particular we give a representation for the universal central extension of $sl (M|N) \otimes A.$

We work with toroidal superalgebra $\tau = T \otimes A \otimes \Omega_{A}/ d_{A}$ with bilinear form $( \ , \ )_T$
on $T$. The algebra $T$ and the form $( \ , \ )_T$ are defined in section $2$.

We will now write down the superbracket on $\tau$ using the superbracket on $T$ given section $2$ $(T1$ to $T10)$\\
$(ST1) \ 1\leq i,j,k,l \leq M$\\ 
$i\neq j, k\neq l$\\
$[T_{ij} (\overline{m}), T_{kl} (\overline{n})]=\begin{cases} 0 \ \mbox{if} \ (\alpha_{ij}, \alpha_{kl}) \geq 0\\
F(\alpha_{ij}, \alpha_{kl}) T_{il} (\overline{m}+\overline{n})\ \mbox{if} \ j=k, l\neq i, (\alpha_{ij}, \alpha_{kl})=-1\\
F(\alpha_{ij}, \alpha_{kl}) T_{kj} (\overline{m}+\overline{n})\ \mbox{if} \ l=i, j\neq k, (\alpha_{ij}, \alpha_{kl})=-1\\
F(\alpha_{ij}, \alpha_{kl}) ((T_{ii} -T_{jj}) (\overline{m}+\overline{n}) + d(t^{\overline{m}})t^{\overline{n}}) \ \mbox{if} \ l=i, j=k, (\alpha_{ij}, \alpha_{kl})=-2\\
\end{cases}$\\[3mm]
$i=j, k\neq l$\\
$[T_{ii}(\overline{m}), T_{kl} (\overline{n})] = (e_{i}, \alpha_{kl}) T_{kl} (\overline{m}+\overline{n}) $\\
$i=j, k=l$\\
$[T_{ii} (\overline{m}) , T_{kk} (\overline{n})] = \delta_{ik}\ d(t^{\overline{m}})t^{\overline{n}}) $\\[3mm]
$$
\begin{array}{llll}
(ST2) & \ 1\leq i, j,k,l \leq N&\\
&[T_{i+M,j+M} (\overline{m}),T_{k+M,l+M} (\overline{n})] =\delta_{jk} \ T_{i+M,l+M} \ (\overline{m}+\overline{n})\\
&-\delta_{li} \ T_{k+M,j+M} \ (\overline{m}+\overline{n})
-\delta_{jk} \ \delta_{li} \ d(t^{\overline{m}})t^{\overline{n}}\\[3mm]
(ST3) & 1\leq i,j,l\leq M,1\leq k \leq N, i \neq j&\\
&[T_{ij} (\overline{m}),T_{k+M,l} (\overline{n})] = \delta_{il} \ F(e_j,e_i) T_{k+M,j} (\overline{m}+\overline{n})\\
&[T_{k+M,l} (\overline{n}),T_{ij} (\overline{m})] = \delta_{il} \ F(e_j,e_i) T_{k+M,j} (\overline{m}+\overline{n})\\
&[T_{ii} (\overline{m}),T_{k+M,l} (\overline{n})] = -\delta_{il} \ T_{k+M,l} (\overline{m}+\overline{n})\\
&[T_{k+M,l} (\overline{n}),T_{ii} (\overline{m})] = \delta_{il} \ T_{k+M,l} (\overline{m}+\overline{n})\\[3mm]
(ST4) & \ 1\leq i, j,k \leq M, 1\leq l \leq N&\\
&[T_{ij} (\overline{m}),T_{k,l+M} (\overline{n})] = \delta_{jk} \ F(e_i,e_j) T_{i,l+M} (\overline{m}+\overline{n})\\
&[T_{k,l+M} (\overline{n}),T_{ij} (\overline{m})] = -\delta_{jk} \ F(e_i,e_j) T_{i,l+M} (\overline{m}+\overline{n})\\[3mm]
(ST5) & \ 1\leq i, j,k \leq N, 1\leq l \leq M&\\
&[T_{i+M,j+M} (\overline{m}),T_{k+M,l} (\overline{n})] =\delta_{jk} \ T_{i+M,l} \ (\overline{m}+\overline{n})\\
&[T_{k+M,l} (\overline{n}),T_{i+M,j+M} (\overline{m})] = -\delta_{jk} \ T_{i+M,l} (\overline{m}+\overline{n})\\[3mm]
(ST6) & \ 1\leq i, j,l \leq N, 1\leq k \leq M&\\
&[T_{i+M,j+M} (\overline{m}),T_{k,l+M} (\overline{n})] =-\delta_{li} \ T_{k,j+M} \ (\overline{m}+\overline{n})\\
&[T_{k,l+M} (\overline{n}),T_{i+M,j+M} (\overline{m})] = \delta_{li} \ T_{k,j+M} (\overline{m}+\overline{n})\\[3mm]
(ST7) & \ 1\leq i,l \leq N, 1\leq j, k\leq M&\\
&[T_{i+M,j} (\overline{m}),T_{k,l+M} (\overline{n})] =\delta_{jk} \ T_{i+M,l+M} (\overline{m}+\overline{n})
+F(e_k,e_j) \delta_{li} \ T_{kj} (\overline{m}+\overline{n})\\
&-\delta_{jk} \ \delta_{li} \ d(t^{\overline{m}})t^{\overline{n}}\\
&[T_{k,l+M} (\overline{n}),T_{i+M,j} (\overline{m})] =\delta_{jk} \ T_{i+M,l+M} (\overline{m}+\overline{n})
+F(e_k,e_j) \delta_{li} \ T_{kj} (\overline{m}+\overline{n})\\
&+\delta_{jk} \ \delta_{li} \ d(t^{\overline{n}})t^{\overline{m}}\\
(ST8) & \ 1\leq i, j\leq M, 1\leq k,l \leq N&\\
&[T_{ij} (\overline{m}),T_{k+M,l+M} (\overline{n})] = [T_{k+M,l+M} (\overline{n}),T_{ij} (\overline{m})] =0\\[3mm]
(ST9) & \ 1\leq i, k\leq N, 1\leq j,l \leq M&\\
&[T_{i+M,j} (\overline{m}),T_{k+M,l} (\overline{n})] = [T_{k+M,l} (\overline{n}),T_{i+M,j} (\overline{m})] =0\\[3mm]
(ST10) & \ 1\leq i, k\leq M, 1\leq j,l \leq N&\\
&[T_{i,j+M} (\overline{m}),T_{k,l+M} (\overline{n})] = [T_{k,l+M} (\overline{n}),T_{i,j+M} (\overline{m})] =0\\[3mm]
\end{array}
$$
\noindent
{\bf (4.3)} We will now extend the vertex operator construction to bigger lattice to accommodate the toroidal case.
Let $\overline{\Gamma}$ be a free $\ZZ$ module on generators $e_1,e_2,\cdots  e_M, \delta_1,\cdots,\delta_{q-1},
d_1,\cdots d_{q-1}.$ (These $\delta_i$'s are not to be confused with $\delta$ function defined earlier. These 
$\delta_i$'s always come with an index).Define a non-degenerate symmetric bilinear form on $\overline{\Gamma}$ by\\
$(e_i,e_j) = \delta_{ij}, (e_i,\delta_j)= (e_i,d_j)=(\delta_i,\delta_j)= (d_i,d_j) = 0,(\delta_i,d_j)= \delta_{ij}$\\
Let $\overline{\hh} = \CC \otimes_{\ZZ} \overline{\Gamma}$ and 
Let $\stackrel{\underline{\wedge}}{\hh} = \oplus_{k \in \ZZ} \overline{\hh} (k) \oplus \CC K$\\
where $\overline{\hh} (k)$ is a copy of $\overline{\hh}$. 
The Lie algebra structure on $\stackrel{\underline{\wedge}}{\hh}$ is given by
$$
[\alpha(k),\beta(m)] = k(\alpha,\beta) \delta_{k+m,0}K
$$

Let $Q$ be the sublattice of $\overline{\Gamma}$ spanned by $e_1,\cdots e_M,\delta_1,\cdots,\delta_{q-1}.$\\
Recall that the lattice $\Gamma$ spanned by $e_1,\cdots e_M$ is defined earlier and hence contained in $\overline{\Gamma}$.\\
Recall that the Co-cycle $F$ is define on $\Gamma$. We now extend the Co-cycle $F$ to $Q$ by
$$
F(e_i,\delta_k) = F(\delta_k,e_i) =F(\delta_k,\delta_l)=1
$$
We further extend $F$ to be a bimultiplicative map $F: Q \times \overline{\Gamma} \rightarrow \{\pm 1 \}
$ is any convenient way.

Consider the group algebra $\CC [\overline{\Gamma}]$ and $\CC [Q]$. Make $\CC [\overline{\Gamma}]$ and $\CC [Q]$ module 
be defining.
$$
e^\alpha . e^\gamma = F(\alpha, \gamma) e^{\alpha +\gamma}, \alpha \in Q, \gamma \in \overline{\Gamma}
$$
Let $\overline{\hh}_{-} = \oplus_{k<0} \ \overline{\hh} (k)$ and consider the Fock space 
$V[\overline{\Gamma}] = \CC [\overline{\Gamma}] \otimes S (\overline{\hh}_-)$ where $S(\overline{\hh}_-)$ is the 
symmetric algebra of $\overline{\hh}_-.$\\
\noindent
{\bf (4.4)} Notation.

Let $\underline{m} = (m_{1},\cdots m_{q-1}) \in \ZZ^{q-1}$\\
So that $\overline{m} = (\underline{m},m_{q}) \in \ZZ^{q}.$

Let $\delta_{\underline{m}} = \displaystyle \sum^{q-1}_{i=1} m_{i} \ \delta_{i}$\\
Define vertex operator as earlier for $\alpha \in Q$
$$
Y(\alpha,z) = e^\alpha z^{\alpha (0)} \ exp \ T_- (\alpha,z) \ exp \ T_+ (\alpha,z)
$$
When $T_+ (\alpha,z)$ are similarly defined as in section $1$. There is a natural $\ZZ_2$-gradation on $\overline{\Gamma}$
comeing from the norm of $\alpha$ in $\overline{\Gamma}$.\\
Write  $\overline{\Gamma} = \overline{\Gamma}_{\overline{0}} \oplus \overline{\Gamma}_{\overline{1}}$\\
For $\alpha \in \overline{\Gamma}_{\overline{0}}$ let $X(\alpha,z) = z^{\frac{(\alpha,\alpha)}{2}} Y(\alpha,z)$\\
For $\alpha \in \overline{\Gamma}_{\overline{1}}$ let $X(\alpha,z) = Y(\alpha,z)$\\
In particular
$$
Y(\delta_{\underline{m}},z)=X(\delta_{\underline{m}},z)=\sum_{n \in \ZZ} X_{n} (\delta_{\underline{m}})z^{-n}
$$
Write \ \ \ \ \ \ \ \ \ \ \ \ \ \ \ \ 
$\alpha (z) =\displaystyle \sum_{n \in \ZZ} \alpha (n) z^{-n-1} , \alpha \in \overline{\Gamma}.$\\[3mm]
The following is very obvious $(\alpha \in Q)$
$$
X (\alpha,z) X (\delta_{\underline{m}},z)=X(\alpha + \delta_{\underline{m}},z)
$$
as \ \ \ \ \ \ \ \ \ \ \ \ \ \ \ \ \ \ \ \ \ \ \ \ \ 
$(\alpha,\delta_{\underline{m}}) = 0$ and $F(\alpha,\delta_{\underline{m}}) = 1.$\\
It is clear that Lemma $1.8$ holds for $\alpha,\beta,\in Q.$

We now consider much bigger Fock space.
$$
V[\overline{\Gamma}] \otimes \FF
$$
The operators $X(\alpha,z),\  \alpha \in Q$ acts on the first component and the operators $\varphi_{i}$'s and
$\varphi^{*}_{i}$ act on the second component. In particular  they commute.\\[3mm]
\noindent
{\bf (4.5)} Recall the $S$ operators and $\alpha_{ij} = e_i - e_j$ from earlier section.\\
Let
\begin{eqnarray*}
S^{\underline{m}}_{ij} (z) =& S_{ij} (z) X (\delta_{\underline{m}}, z)\\
=&\displaystyle \sum_{n \in \ZZ} S^{\underline{m}}_{ij} (n) z^{-n-1}
\end{eqnarray*}
so that 
$S^{\underline{m}}_{ij} (n)= \displaystyle \sum_{k \in \ZZ} S_{ij} (k) X_{n-k} (\delta_{\underline{m}})$\\[2mm]
Let
$$
\alpha (z) X (\delta_{\underline{m}} ,z)= \displaystyle \sum_{n \in \ZZ} T^{\alpha}_{n} (\delta_{\underline{m}}) z^{-n-1}
$$
so that 
$ T^{\alpha}_{n} (\delta_{\underline{m}})= \displaystyle \sum_{k \in \ZZ} \alpha(k) X_{n-k} (\delta_{\underline{m}})$\\

We now state and Prove the main theorem of this paper
\paragraph*{Theorem (4.6)}
The operators $X$ and $S$ define a representation of \\
$\tau = T\otimes A \otimes \Omega_A/d_n$ on the space $V(\overline{\Gamma}) \otimes \FF$ by the following correspondence\\
$1\leq i,j\leq M$\\
$i\neq j \ \ \ \ \ \ \ \ \ T_{ij} \otimes t^{\overline{m}} \mapsto X_{m_q} (\alpha_{ij} +\delta_{\underline{m}})$\\
$i= j \ \ \ \ \ \ \ \ \ T_{ii} \otimes t^{\overline{m}} \mapsto T^{e_{i}}_{m_q} (\delta_{\underline{m}})$\\
If one of $i$ or $j>M$\\
$T_{ij} \otimes t^{\overline{m}} \mapsto S^{\underline{m}}_{ij} (m_q)$\\
$d(t^{\overline{m}})t^{\overline{n}} \mapsto X_{{m_q}+{n_q}} (\delta_{\underline{m} + \underline{n}}) + {m_q} \ 
T^{\delta_{\underline{m}}}_{{m_q}+{n_q}} (\delta_{\underline{m} +\underline{n}})$\\
We need the following
\paragraph*{Lemma 4.7}
$$
\begin{array}{llll}
(1)& D_{z} z X (\delta_{\underline{m}},z) &=& z \delta_{\underline{m}} (z) X (\delta_{\underline{m}},z) + X (\delta_{\underline{m}},z)\\
(2)& D_{z}  X (\delta_{\underline{m}},z) &=& \delta_{\underline{m}} (z) X (\delta_{\underline{m}},z) \\
\end{array}
$$
\paragraph*{Proof} Follows from Lemma(4.6) of $[EM]$ just note that $\alpha(z)$ in our paper and $\alpha(z)$  in $[EM]$
differ by a factor of $z$.
\paragraph*{Lemma 4.8}
$$
\begin{array}{llll}
(1)&D_z \delta ( z/w). zX (\delta_{\underline{m}},z) X (\delta_{\underline{n}},w)\\
&= wX (\delta_{\underline{m} + \underline{n}},w) D_z \ \delta (z/w) \ \ \ \ \\
&- w \ \delta_{\underline{m}} (w) X(\delta_{{\underline{m}} +{\underline{n}}},w) \delta(z/w)\\
&- X(\delta_{{\underline{m}} +{\underline{n}}},w) \delta(z/w)\\[2mm]
(2)&D_z \delta ( z/w) X (\delta_{\underline{m}},z) X (\delta_{\underline{n}},w)\\
&= X(\delta_{{\underline{m}} +{\underline{n}}},w) D_{z} \delta(z/w)\\
&- \delta_{\underline{m}} (w) X(\delta_{{\underline{m}} +{\underline{n}}},w) \delta(z/w)\\
\end{array}
$$
\paragraph*{Proof} Follows from Lemma $1.5(3)$ and Lemma $(4.7)$
\paragraph*{Lemma 4.9}
$$
T^{\delta_{\underline{m}}}_{m_q} (\delta_{\underline{m}}) + {m_q} \ X_{m_q} (\delta_{\underline{m}}) =0
$$
Follows from Lemma 3.13 of $[EM].$\\
\noindent
{\bf (4.10)} We have noted that $d(t^{\overline{m}})t^{\overline{n}} = -d(t^{\overline{n}})t^{\overline{m}}.$
Because of the Lemma $4.9$ this can be verified in the statement of the Theorem $(4.6)$.\\
Proof of the Theorem $(4.6).$

It is sufficient to verify the following relation among the $X$ and $S$ operators.
$$
\begin{array}{llll}
(SR1) & 1\leq i, j,k,l \leq M&\\
(1) & i\neq j, k\neq l\\
&[X(\alpha_{ij},z) X(\delta_{\underline{m}},z),X(\alpha_{kl},w) X(\delta_{\underline{n}},w]
=o\  \mbox{if} \ (\alpha_{ij},\alpha_{kl}) \geq 0\\
&=F(\alpha_{ij},\alpha_{kl}) X (\alpha_{ij}+\alpha_{kl},z) X(\delta_{{\underline{m}} +{\underline{n}}},z) \delta (z/w)\  \mbox{if} \ (\alpha_{ij},\alpha_{kl}) =-1\\
&=F(\alpha_{ij},-\alpha_{ij})\alpha_{ij} (w) w X(\delta_{{\underline{m}} +{\underline{n}}},w) \delta (z/w)\\
&-F(\alpha_{ij},-\alpha_{ij}) (wX(\delta_{{\underline{m}} +{\underline{n}}},w) \ D_{z} \delta (z/w)\\
&-w \ \delta_{\underline{m}} (w) X(\delta_{{\underline{m}} +{\underline{n}}},w) \delta (z/w)\\
&-X(\delta_{{\underline{m}} +{\underline{n}}},w) \delta (z/w)) \ \mbox{if} \ \ (\alpha_{ij},\alpha_{kl}) =-2\\[3mm]
(2) & i = j, k \neq l\\
&[e_i(z) X(\delta_{\underline{m}},z),X(\alpha_{kl},w) X(\delta_{\underline{n}},w]\\
&= (e_i,\alpha_{kl}) w^{-1} X(\alpha_{kl},w) X(\delta_{{\underline{m}} +{\underline{n}}},w) \delta (z/w)\\[3mm]
(3) & i = j, k = l\\
&=[e_i(z) X(\delta_{\underline{m}},z), e_k (w) X (\delta_{\underline{n}},w)]\\
&= -\delta_{kl} \ w^{-1}  X(\delta_{{\underline{m}} +{\underline{n}}},w) D_z \ \delta (z/w)\\
& +\delta_{kl} \ w^{-1} \ \delta_{\underline{m}} (w) X(\delta_{{\underline{m}} +{\underline{n}}},w) \delta (z/w)\\[3mm]
(SR2) & 1\leq i, j,k,l \leq N&\\
&[S^{\underline{m}}_{i+M,j+M} (z), S^{\underline{n}}_{k+M,l+M} (w)]\\
&=\delta_{jk} \ S^{{\underline{m}}+{\underline{n}}}_{i+M,l+M} (w) w^{-1} \ \delta (z/w)\\
&-\delta_{li} \ S^{{\underline{m}}+{\underline{n}}}_{k+M,l+M} (w) w^{-1} \ \delta (z/w)\\
&+\delta_{jk} \delta_{li} (w^{-1} X(\delta_{{\underline{m}} +{\underline{n}}},w) \ D_{z} \delta (z/w)\\
&- w^{-1} \delta_{\underline{m}} (w) X(\delta_{{\underline{m}} +{\underline{n}}},w)\  \delta (z/w))\\[3mm]
\end{array}
$$
\newpage
$$
\begin{array}{llll}
(SR3) & 1\leq i,l \leq M,1\leq k \leq N&\\
(1) & i\neq j\\
&[X(\alpha_{ij},z) X(\delta_{\underline{m}},z) , S^{\underline{n}}_{k+M,l} (w)]\\
&=\delta_{il} \ F(e_j,e_i) S^{{\underline{m}}+{\underline{n}}}_{k+M,j} (z) \ \delta (z/w)\\
&[S^{\underline{n}}_{k+M,l} (z), X(\alpha_{ij},w) X(\delta_{\underline{m}},w)]\\
&=\delta_{il} \ F(e_i,e_j) S^{{\underline{m}}+{\underline{n}}}_{k+M,l} (w) \ \delta (z/w)\\[3mm]
(2) & i= j\\
&[e_i(z) X(\delta_{\underline{m}},z),S^{\underline{n}}_{k+M,l} (w)]\\
&=-\delta_{il} \ S^{{\underline{m}}+{\underline{n}}}_{k+M,l} (w) w^{-1} \ \delta (z/w)\\
&[S^{\underline{n}}_{k+M,l} (z),e_{i}(w) X(\delta_{\underline{m}},w)]\\
&=\delta_{il} \ S^{{\underline{m}}+{\underline{n}}}_{k+M,l} (w) w^{-1} \ \delta (z/w)\\[3mm]
(SR4) & 1\leq i,j,k \leq M,1\leq l \leq N&\\
(1) & i\neq j\\
&[X(\alpha_{ij},z) X(\delta_{\underline{m}},z) , S^{\underline{n}}_{k,l+M} (w)]\\
&=\delta_{jk} \ F(e_i,e_j) S^{{\underline{m}}+{\underline{n}}}_{i,l+M} (w) \ \delta (z/w)\\
&[S^{\underline{n}}_{k,l+M} (z), X(\alpha_{ij},w) X(\delta_{\underline{m}},w)]\\
&=-\delta_{jk} \ F(e_i,e_j) S^{{\underline{m}}+{\underline{n}}}_{i,l+M} (w) \ \delta (z/w)\\[3mm]
(2) & i= j\\
&[e_i(z) X(\delta_{\underline{m}},z),S^{\underline{n}}_{k,l+M} (w)]\\
&=\delta_{ik} \ S^{{\underline{m}}+{\underline{n}}}_{i,l+M} (w) w^{-1} \ \delta (z/w)\\
&[S^{\underline{n}}_{k,l+M} (z),e_{i}(w) X(\delta_{\underline{m}},w)]\\
&=-\delta_{ik} \ S^{{\underline{m}}+{\underline{n}}}_{i,l+M} \  w^{-1} \ \delta (z/w)\\[3mm]
(SR5) & 1\leq i,j,k \leq N,1\leq l \leq M&\\
&[S^{\underline{m}}_{i+M,j+M} (z), S^{\underline{n}}_{k+M,l} (w)]\\
&=\delta_{jk} \ S^{{\underline{m}}+{\underline{n}}}_{i+M,l} (w) w^{-1} \ \delta (z/w)\\
&[S^{\underline{n}}_{k+M,l} (z), S^{\underline{m}}_{i+M,j+M} (w)]\\
&=-\delta_{jk} \ S^{{\underline{m}}+{\underline{n}}}_{i+M,l} \  w^{-1} \ \delta (z/w)\\[3mm]
\end{array}
$$
\newpage
$$
\begin{array}{llll}
(SR6) & 1\leq i, j,l \leq N, 1\leq k \leq M&\\
&[S^{\underline{m}}_{i+M,j+M} (z), S^{\underline{n}}_{k,l+M} (w)] 
=-\delta_{li} \ S^{{\underline{m}}+{\underline{n}}}_{k,j+M} (w) w^{-1} \ \delta (z/w)\\
&[S^{\underline{n}}_{k,l+M} (z), S^{\underline{m}}_{i+M,j+M} (w)] 
=\delta_{li} \ S^{{\underline{m}}+{\underline{n}}}_{k,j+M} (w) w^{-1} \ \delta (z/w)\\[3mm]
(SR7) & 1\leq i,l \leq N,1\leq j,k \leq M&\\
& k\neq j\\
&[S^{\underline{m}}_{i+M,j} (z), S^{\underline{n}}_{k,l+M} (w)]\\
&=\delta_{jk} \ S^{{\underline{m}}+{\underline{n}}}_{i+M,l+M} (w)w^{-1} \ \delta (z/w)\\
&+\delta_{li} \ F(e_k,e_j) X(\alpha_{kj} + \delta_{{\underline{m}} +{\underline{n}}},w) w^{-2} \ \delta (z/w)\\
&+\delta_{jk} \delta_{li} (w^{-1} X(\delta_{{\underline{m}} +{\underline{n}}},w) \ D_{z} \delta (z/w)\\
&- w^{-1} \delta_{\underline{m}} (w) X(\delta_{{\underline{m}} +{\underline{n}}},w)\  \delta (z/w))\\[3mm]
& k\neq j\\
&[S^{\underline{n}}_{k,l+M} (z), S^{\underline{m}}_{i+M,j} (w)]\\
&=\delta_{jk} \ S^{{\underline{m}}+{\underline{n}}}_{i+M,l+M} w^{-1} \ \delta (z/w)\\
&+\delta_{li} \ F(e_k,e_j) X(\alpha_{kj} + \delta_{{\underline{m}} +{\underline{n}}},w) w^{-2} \ \delta (z/w)\\
&-\delta_{jk} \delta_{li} (w^{-1} X(\delta_{{\underline{m}} +{\underline{n}}},w) \ D_{z} \delta (z/w)\\
&- w^{-1} \delta_{\underline{n}} (w) X(\delta_{{\underline{m}} +{\underline{n}}},w)\  \delta (z/w))\\[3mm]
& j=k\\
&[S^{\underline{m}}_{i+M,j} (z), S^{\underline{n}}_{k,l+M} (w)]\\
&=\delta_{jk} \ S^{{\underline{m}}+{\underline{n}}}_{i+M,l+M} (w) w^{-1} \ \delta (z/w)\\
&+\delta_{li} \ F(e_k,e_j) e_{j} (w) w^{-1} X (\delta_{{\underline{m}} +{\underline{n}}},w) \ \delta (z/w)\\
&+\delta_{jk} \ \delta_{li} \ w^{-1} X(\delta_{{\underline{m}} +{\underline{n}}},w) \ D_{z} \delta (z/w)\\
&- \delta_{jk} \ \delta_{li} \ w^{-1} \  \delta_{\underline{m}} (w) X(\delta_{{\underline{m}} +{\underline{n}}},w)\ 
\delta (z/w)\\[3mm]
& j=k\\
&[S^{\underline{n}}_{k,l+M} (z), S^{\underline{m}}_{i+M,j} (w)]\\
&=\delta_{jk} \ S^{{\underline{m}}+{\underline{n}}}_{i+M,l+M} (w) w^{-1} \ \delta (z/w)\\
&+\delta_{li} \ F(e_k,e_j) e_{j} (w) w^{-1} X (\delta_{{\underline{m}} +{\underline{n}}},w) \ \delta (z/w)\\
&-\delta_{jk} \ \delta_{li} \ w^{-1} X(\delta_{{\underline{m}} +{\underline{n}}},w) \ D_{z} \delta (z/w)\\
&+ \delta_{jk} \ \delta_{li} \ w^{-1} \  \delta_{\underline{n}} (w) X(\delta_{{\underline{m}} +{\underline{n}}},w)\ \delta (z/w)\\[3mm]\\
\end{array}
$$
\newpage
$$
\begin{array}{llll}
(SR8) & 1\leq i,j \leq M,1\leq k,l \leq N&\\
&[X(\alpha_{ij},z) X(\delta_{\underline{m}},z) , S^{\underline{n}}_{k+M,l+M} (w)] =\\
&[S^{\underline{n}}_{k+M,l+M} (z), X(\alpha_{ij},w) X(\delta_{\underline{m}},w)] = 0\\
\end{array}
$$
Note that for $i=j,\ X(\alpha_{ij},z)$ has to be replaced by $e_i(z)$.\\
$$
\begin{array}{llll}
(SR9) & 1\leq i,k \leq N,1\leq j,l \leq M&\\
&[S^{\underline{m}}_{i+M,j} (z), S^{\underline{n}}_{k+M,l} (w)]=\\
&[S^{\underline{m}}_{k+M,l} (z), S^{\underline{n}}_{i+M,j} (w)]=0\\[3mm]
(SR10) & 1\leq i,k \leq M,1\leq j,l \leq N&\\
&[S^{\underline{m}}_{i,j+M} (z), S^{\underline{n}}_{k,l+M} (w)]\\
&=[S^{\underline{n}}_{k,l+M} (z), S^{\underline{m}}_{i,j+M} (w)]=0\\[3mm]
\end{array}
$$
When we expand the above infinite series in components form we will see that the relations $SR1$ to $SR10$
will be equal to $ST1$ to $ST10$. In order to verify the relations $SR1$ to $SR10$, we need to rewrite the
relations $R1$ to $R10$ in the infinite series form. Thus we have the following.\\
$$
\begin{array}{llll}
(S1) & 1\leq i, j,k,l \leq M&\\
(1) & i\neq j, k\neq l\\
&[X(\alpha_{ij},z) ,X(\alpha_{kl},w)] = 0 \ \mbox{if} \ (\alpha_{ij},\alpha_{kl}) \geq 0\\
&=F(\alpha_{ij},\alpha_{kl}) X (\alpha_{ij}+\alpha_{kl},z) \delta (z/w)\  \mbox{if} \ (\alpha_{ij},\alpha_{kl}) =-1\\
&=F(\alpha_{ij},-\alpha_{ij}) (\alpha_{ij} (z) z \ \delta (z/w) -
z\ D_{z} \delta (z/w)) \ \mbox{if} \ \ (\alpha_{ij},\alpha_{kl}) =-2\\[3mm]
(2) & i = j, k \neq l\\
&[e_i(z),X(\alpha_{kl},w)] 
= z^{-1} (e_i,\alpha_{kl}) X(\alpha_{kl},z) \ \delta (z/w)\\[3mm]
(3) & i = j, k = l\\
&[e_i(z) , e_k (w)] = -\delta_{ik} \ w^{-1} \  D_z \ \delta (z/w)\\[3mm]
\end{array}
$$
\newpage
$$
\begin{array}{llll}
(S2) & 1\leq i, j,k,l \leq N&\\
&[S_{i+M,j+M} (z), S_{k+M,l+M} (w)]=\\
&\delta_{jk} \ S_{i+M,l+M} (z)\  w^{-1} \ \delta (z/w)\\
&-\delta_{li} \ S_{k+M,j+M} (z)\  w^{-1} \ \delta (z/w)\\
&+\delta_{jk} \ \delta_{li} \ w^{-1} \ \ D_{z} \delta (z/w)\\[3mm]
(S3) & 1\leq i,j,l \leq M,1\leq k \leq N&\\
(1) & i\neq j\\
&[X(\alpha_{ij},z), S_{k+M,l} (w)]\\
&=\delta_{il} \ F(e_j,e_i) S_{k+M,j} (z) \ \delta (z/w)\\
&[S_{k+M,l} (z), X(\alpha_{ij},w)]\\
&= \delta_{il} \ F(e_i,e_j) S_{k+M,j} (z) \ \delta (z/w)\\[3mm]
(2) & i= j\\
&[e_i(z),S_{k+M,l} (w)]\\
&=-\delta_{il} \ S_{k+M,l} (z)\  w^{-1} \ \delta (z/w)\\
&[S_{k+M,l} (z),e_{i}(z)]\\
&=\delta_{il} \ S_{k+M,l} (z) \ w^{-1} \ \delta (z/w)\\[3mm]
(S4) & 1\leq i,j,k \leq M,1\leq l \leq N&\\
(1) & i\neq j\\
&[X(\alpha_{ij},z), S_{k,l+M} (w)]\\
&=\delta_{jk} \ F(e_i,e_j) S_{i,l+M} (z) \ \delta (z/w)\\
&[S_{k,l+M} (z), X(\alpha_{ij},z)]\\
&=-\delta_{jk} \ F(e_i,e_j) S_{i,l+M} \ (z)\ \delta (z/w)\\[3mm]
(2) & i= j\\
&[e_i(z),S_{k,l+M} (w)]\\
&=\delta_{ik} \ S_{i,l+M} (z)\  w^{-1} \ \delta (z/w)\\
&[S_{k,l+M} (z),e_{i}(w)]\\
&=-\delta_{ik} \ S_{i,l+M} (z)\  w^{-1} \ \delta (z/w)\\[3mm]
\end{array}
$$
\newpage
$$
\begin{array}{llll}
(S5) & 1\leq i,j,k \leq N,1\leq l \leq M&\\
&[S_{i+M,j+M} (z), S_{k+M,l} (w)]\\
&=\delta_{jk} \ S_{i+M,l} (z)\  w^{-1} \ \delta (z/w)\\
&[S_{k+M,l} (z), S_{i+M,j+M} (w)]\\
&=-\delta_{jk} \ S_{i+M,l} (z)\  w^{-1} \ \delta (z/w)\\[3mm]
(S6) & 1\leq i, j,l \leq N, 1\leq k \leq M&\\
&[S_{i+M,j+M} (z), S_{k,l+M} (w)]\\ 
&=-\delta_{li} \ S_{k,j+M} (z)\  w^{-1} \ \delta (z/w)\\
&[S_{k,l+M} (z), S_{i+M,j+M} (w)]\\
&=\delta_{li} \ S_{k,j+M} (z)\  w^{-1} \ \delta (z/w)\\[3mm]
(S7) & 1\leq i,l \leq N,1\leq j,k \leq M&\\
(1) & j\neq k\\
&[S_{i+M,j} (z), S_{k,l+M} (w)]\\
&=\delta_{jk} \ S_{i+M,l+M} (z)\  w^{-1} \ \delta (z/w)\\
&+\delta_{li} \ F(e_k,e_j) X(\alpha_{kj},z) z^{-1} w^{-1} \ \delta (z/w)\\
&+\delta_{kj} \ \delta_{li} \ w^{-1}  \ D_{z} \delta (z/w)\\[3mm]
& j = k\\
&[S_{i+M,j} (z), S_{k,l+M} (w)]\\
&=\delta_{jk} \ S_{i+M,l+M} (z) w^{-1} \ \delta (z/w)\\
&+\delta_{li} \ F(e_k,e_j) e_{j} \ (z) \ w^{-1} \ \delta (z/w)\\
&+\delta_{kj} \ \delta_{li} \ w^{-1}  \ D_{z} \delta (z/w)\\[3mm]
(2)& j \neq k\\
&[S_{k,l+M} (z), S_{i+M,j} (w)]\\
&=\delta_{jk} \ S_{i+M,l+M} (z) w^{-1} \ \delta (z/w)\\
&+\delta_{li} \ F(e_k,e_j) X(\alpha_{kj},z) z^{-1} w^{-1}  \ \delta (z/w)\\
&-\delta_{jk} \ \delta_{li} \ w^{-1}  \ D_{z} \delta (z/w)\\[3mm]
\end{array}
$$
\newpage
$$
\begin{array}{llll}
& j=k\\
&[S_{k,l+M} (z), S_{i+M,j} (w)]\\
&=\delta_{jk} \ S_{i+M,l+M} (z) \ w^{-1} \ \delta (z/w)\\
&+\delta_{li} \ F(e_k,e_j) e_{j} (z) w^{-1} \ \delta (z/w)\\
&-\delta_{jk} \ \delta_{li} \ w^{-1} \ D_{z} \delta (z/w)\\[3mm]
(S8) & 1\leq i,j \leq M,1\leq k, l \leq N&\\
&[X(\alpha_{ij},z) , S_{k+M,l+M} (w)] =
[S_{k+M,l+M} (z), X(\alpha_{ij},w)] = 0\\
\end{array}
$$
Note that for $i=j,\ X(\alpha_{ij},z)$ has to be replaced by $e_i(z)$.\\
$$
\begin{array}{llll}
(S9) & 1\leq i,k \leq N,1\leq j,l \leq M&\\
&[S_{i+M,j} (z), S_{k+M,l} (w)]=
[S_{k+M,l} (z), S_{i+M,j} (w)]=0\\[3mm]
(S10) & 1\leq i,k \leq M,1\leq j,l \leq N&\\
&[S_{i,j+M} (z), S_{k,l+M} (w)]
=[S_{k,l+M} (z), S_{i,j+M} (w)]=0\\[3mm]
\end{array}
$$
It is now a simple matter to verify $SR1$ to $SR10$ by using Lemma $4.7$ and $4.8$ and the relations $S1$ to $S10$.
We will only verify the third part of $SR1$ and $SR7$ which is more complex then others.
The rest of the relations can be verified similarly.

Recall that the operators $X(\delta_{\underline{m}},z)$ commute with all operators. To verify $SR1$ (third part)
consider\\
$$
\begin{array}{llll}
&[e_i(z) X(\delta_{\underline{m}},z),e_k(w) X(\delta_{\underline{n}},w)]\\
&= [e_i (z),e_k (w)] X (\delta_{\underline{m}},z) X(\delta_{\underline{n}},w)\\
&= -\delta_{ki} \ w^{-1} D_z \ \delta (z/w) X(\delta_{\underline{m}} ,z) X(\delta_{\underline{n}},w) \ (by \ S1 \ )\\
&= -\delta_{ki} \ w^{-1}  X(\delta_{{\underline{m}} +{\underline{n}}},w) D_z \ \delta (z/w)\\
& +\delta_{ki} \ w^{-1} \ \delta_{\underline{m}} (w) X(\delta_{{\underline{m}} +{\underline{n}}},w) \ \delta (z/w) \ (by \ Lemma \ 4.8(2) \ )\\[3mm]
\end{array}
$$
Which is precisely $SR1$ (3rd part)\\
\newpage
We will now verify the first and third part of $SR7$. Consider for $k\neq j$\\
$$
\begin{array}{llll}
&[S^{\underline{m}}_{i+M,j} (z), S^{\underline{n}}_{k,l+M} (w)]\\
&=[S_{i+M,j} (z) X(\delta_{\underline{m}} ,z),S_{k,l+M} (w) X(\delta_{\underline{n}},w)] \ (by \ 4.5 \ )\\
&=[S_{i+M,j} (z),S_{k,l+M} (w)] \ X(\delta_{\underline{m}} ,z)\ X(\delta_{\underline{n}},w)\\
&=\delta_{jk} \ S_{i+M,l+M} (z) w^{-1} \ \delta (z/w) \ X(\delta_{\underline{m}} ,z)\ X(\delta_{\underline{n}},w)\\
&+\delta_{li} \ F(e_k,e_j) X(\alpha_{kj},z) z^{-1} w^{-1}  \ \delta (z/w) X(\delta_{\underline{m}} ,z)\ X(\delta_{\underline{n}},w)\\
&+\delta_{jk} \ \delta_{li} \ w^{-1}  \ D_{z} \delta (z/w) X(\delta_{\underline{m}} ,z)\ X(\delta_{\underline{n}},w)\ (by \ S7 \ ) \\
&=\delta_{jk} \ S_{i+M,l+M} (w)\ w^{-1} X(\delta_{\underline{m}} ,w)\ X(\delta_{\underline{n}},w)\ \delta (z/w) \\
&+\delta_{li} \ F(e_k,e_j) X(\alpha_{kj},w) w^{-2} \  X(\delta_{\underline{m}} ,w)\ X(\delta_{\underline{n}},w) \ \delta (z/w)\\
&+\delta_{jk} \delta_{li} (w^{-1} X(\delta_{{\underline{m}} +{\underline{n}}},w) \ D_{z} \delta (z/w)\\
&- w^{-1} \delta_{\underline{m}} (w) X(\delta_{{\underline{m}} +{\underline{n}}},w)\  \delta (z/w))
\ \ (By \ Lemma \ 4.7 \ and \ 4.8 \ )\\[3mm]
\end{array}
$$
for $k = j$\\
$$
\begin{array}{llll}
&[S^{\underline{m}}_{i+M,j} (z), S^{\underline{n}}_{k,l+M} (w)]
=[S_{i+M,j} (z) X(\delta_{\underline{m}} ,z),S_{k,l+M} (w) X(\delta_{\underline{n}},w)] \ (by \ 4.5 \ )\\
&=[S_{i+M,j} (z),S_{k,l+M} (w)] \ X(\delta_{\underline{m}} ,z)\ X(\delta_{\underline{n}},w)\\
&=\delta_{jk} \ S_{i+M,l+M} (z) w^{-1} \ \delta (z/w) \ X(\delta_{\underline{m}} ,z)\ X(\delta_{\underline{n}},w)\\
&+\delta_{li} \ F(e_k,e_j) e_j (z) w^{-1}  \ \delta (z/w) X(\delta_{\underline{m}} ,z)\ X(\delta_{\underline{n}},w)\\
&+\delta_{jk} \ \delta_{li} \ w^{-1}  \ D_{z} \delta (z/w) X(\delta_{\underline{m}} ,z)\ X(\delta_{\underline{n}},w)\ (by \ S7 \ ) \\
&=\delta_{jk} \ S_{i+M,l+M} (w)\ w^{-1} X(\delta_{\underline{m}}+\delta_{\underline{n}},w)\ \delta (z/w) \\
&+\delta_{li} \ F(e_k,e_j) e_j \ (w) w^{-1} \  X(\delta_{\underline{m}} +\delta_{\underline{n}},w) \ \delta (z/w)\\
&+\delta_{jk} \ \delta_{li} \ w^{-1} X(\delta_{{\underline{m}} +{\underline{n}}},w) \ D_{z} \delta (z/w)\\
&-\delta_{jk} \ \delta_{li} \  w^{-1} \delta_{\underline{m}} (w) X(\delta_{{\underline{m}} +{\underline{n}}},w)\  \delta (z/w)
\ \ (By \ Lemma \ 4.7 \ and \ 4.8 \ )\\[3mm]
\end{array}
$$
\newpage
We now verify the 2nd part and 4th part of $SR7$ Consider for $k\neq j$\\
$$
\begin{array}{llll}
&[S^{\underline{n}}_{k,l+M} (z), S^{\underline{m}}_{i+M,j} (w)]
=[S_{k,l+M} (z) X(\delta_{\underline{n}} ,z),S_{i+M,j} (w) X(\delta_{\underline{m}},w)] \ (by \ 4.5 \ )\\
&=[S_{k,l+M} (z),S_{i+M,j} (w)] \ X(\delta_{\underline{n}} ,z)\ X(\delta_{\underline{m}},w)\\
&=\delta_{jk} \ S_{i+M,l+M} (z) w^{-1} \ \delta (z/w) \ X(\delta_{\underline{n}} ,z)\ X(\delta_{\underline{m}},w)\\
&+\delta_{li} \ F(e_k,e_j) X(\alpha_{kj},z)\  z^{-1} \  w^{-1}  \ \delta (z/w) X(\delta_{\underline{n}} ,z)\ X(\delta_{\underline{m}},w)\\
&-\delta_{jk} \ \delta_{li} \ w^{-1}  \ D_{z} \delta (z/w) X(\delta_{\underline{n}} ,z)\ X(\delta_{\underline{m}},w)\ (by \ S7 \ ) \\
&=\delta_{jk} \ S_{i+M,l+M} (w)\ w^{-1} X(\delta_{\underline{m}}+\delta_{\underline{n}},w)\ \delta (z/w) \\
&+\delta_{li} \ F(e_k,e_j) \ X(\alpha_{kj},w) \ w^{-2} \ \delta (z/w) X(\delta_{\underline{m}} +\delta_{\underline{n}},w) \\
&-\delta_{jk} \ \delta_{li} \  X(\delta_{{\underline{m}} +{\underline{n}}},w) \ w^{-1}\ D_{z} \delta (z/w)\\
&+\delta_{jk} \ \delta_{li} \  w^{-1} \delta_{\underline{n}} (w) X(\delta_{{\underline{m}} +{\underline{n}}},w)\  \delta (z/w))\\[3mm]
\end{array}
$$
for $k = j$\\
$$
\begin{array}{llll}
&[S^{\underline{n}}_{k,l+M} (z), S^{\underline{m}}_{i+M,j} (w)]
=[S_{k,l+M} (z) X(\delta_{\underline{n}} ,z),S_{i+M,j} (w) X(\delta_{\underline{m}},w)] \ (by \ 4.5 \ )\\
&=[S_{k,l+M,j} (z),S_{i+M,j} (w)] \ X(\delta_{\underline{n}} ,z)\ X(\delta_{\underline{m}},w)\\
&=\delta_{jk} \ S_{i+M,l+M} (z) w^{-1} \ \delta (z/w) \ X(\delta_{\underline{n}} ,z)\ X(\delta_{\underline{m}},w)\\
&+\delta_{li} \ F(e_k,e_j) \ e_{j} (z)\  w^{-1}  \ \delta (z/w) X(\delta_{\underline{n}} ,z)\ X(\delta_{\underline{m}},w)\\
&-\delta_{jk} \ \delta_{li} \ w^{-1}  \ D_{z} \delta (z/w) X(\delta_{\underline{n}} ,z)\ X(\delta_{\underline{m}},w)\ (by \ S7 \ ) \\
&=\delta_{jk} \ S_{i+M,l+M} (z)\ w^{-1} X(\delta_{\underline{m}}+\delta_{\underline{n}},w)\ \delta (z/w) \\
&+\delta_{li} \ F(e_k,e_j) \ e_{j} (w) \ w^{-1} \ X(\delta_{\underline{m}} +\delta_{\underline{n}},w) \ \delta (z/w) \\
&-\delta_{jk} \ \delta_{li} \  w^{-1} \ X(\delta_{{\underline{m}} +{\underline{n}}},w) \ D_{z} \delta (z/w)\\
&+\delta_{jk} \ \delta_{li} \  w^{-1} \delta_{\underline{n}} (w) X(\delta_{{\underline{m}} +{\underline{n}}},w)\  \delta (z/w))\\[3mm]
\end{array}
$$
Note that we are useing Lemmas $4.7$ and $4.8$ in the above calculations.
\newpage
\paragraph{Remark 4.10}
\begin{enumerate}
\item  The central element $t^{\overline{m}} K_i, 1\leq i \leq q-1$ acts as $T^{\delta_{i}}_{{m_{q}}} (\delta_{\underline{m}}).$
The central element  $t^{\overline{m}} K_q$ act as $X_{m_q} (\delta_{\underline{m}}).$\ In particular $K_q$ acts as $X_{0} (0) = I_d.$
\item The module $V[\overline{\Gamma}] \otimes \FF$ \ is integrable with respect to $\stackrel{\wedge}{sl(M)}$ but not
integrable with respect to  $\stackrel{\wedge}{sl(N)}$ the second component of the even part. See $[EZ]$ and $[EF]$ for some interesting results
on integrable modules.
\item The module $V[\overline{\Gamma}] \otimes \FF$ \  is not irreducible.\\
 For example $\delta_{\underline{m}} (k) (V[\overline{\Gamma}] \otimes \FF)$ is always a proper sub module for any \\
$\underline{m} \in \ZZ^{q-1}$ and any $k < 0.$
\item The structure of submodules and quotient modules of $V[\overline{\Gamma}] \otimes \FF$ \  will be investigated in a 
subsequent paper.
\item Theorem (4.5), when restricted to $\stackrel{\wedge}{sl(M)}$, recover the main result of $[EM]$ and $[EMY]$.
  In fact the proof here as much simpler than $[EM]$. In the Lie algebra case our proof works for any simply laced case.
\item It is easy to see the representation is faithful. The only nontrivial thing is to prove that the full center acts
faithfully. This follows from $[EM]$.\\[2mm]
\end{enumerate}
\emph{Acknowledgements}: I thank Drazen Adamovic for some clarification on super vertex operator algebras.

\end{document}